\documentclass{siamart251216}

\usepackage{amsfonts}
\usepackage{graphicx}
\usepackage{epstopdf}
\usepackage{algorithmic}
\usepackage{bbm}
\usepackage{tikz}
\usepackage{lipsum}
\usepackage{multirow}
\usetikzlibrary{arrows}
\usepackage{enumitem}
\ifpdf
  \DeclareGraphicsExtensions{.eps,.pdf,.png,.jpg}
\else
  \DeclareGraphicsExtensions{.eps}
\fi

\makeatletter
\def\namedlabel#1#2{\begingroup
    #2%
    \def\@currentlabel{#2}%
    \phantomsection\label{#1}\endgroup
}
\makeatother

\usepackage{cleveref}
\crefname{proposition}{Proposition}{Propositions}
\Crefname{proposition}{Proposition}{Propositions}

\newsiamremark{remark}{Remark}
\newsiamremark{hypothesis}{Hypothesis}
\crefname{hypothesis}{Hypothesis}{Hypotheses}
\newsiamthm{claim}{Claim}
\newsiamthm{assumption}{Assumption}

\headers{}{}
\title{Policy stability and ultimate stationarity in discounted risk-sensitive stochastic control}

 \author{
Nicole B{\"a}uerle\thanks{Department of Mathematics, Karlsruhe Institute of Technology (KIT), Karlsruhe, Germany (\email{nicole.baeuerle@kit.edu}).}
 \and
Marcin Pitera\thanks{Institute of Mathematics, Jagiellonian University, Krakow, Poland
  (\email{marcin.pitera@uj.edu.pl}).}
  \and
{\L}ukasz Stettner\thanks{Institute of Mathematics, Polish Academy of Sciences, Warsaw, Poland
  (\email{l.stettner@impan.pl}).}}

\usepackage{amsopn}
\usepackage{enumitem}

\DeclareMathOperator*{\esssup}{ess\,sup} % ess sup
\DeclareMathOperator*{\essinf}{ess\,inf} % ess inf

\newcommand{\1}{\mathbbm{1}}     

\DeclareMathOperator*{\argmax}{arg\,max}
\def\cF{\mathcal{F}}

\def\bP{\mathbb{P}}
\def\bE{\mathbb{E}}
\def\bR{\mathbb{R}}

\def\bN{\mathbb{N}}

\def\bC{\mathbb{C}}

\def\cP{\mathcal{P}}

\makeatletter
\def\namedlabel#1#2{\begingroup
    #2%
    \def\@currentlabel{#2}%
    \phantomsection\label{#1}\endgroup
}
\makeatother

\newtheorem{example}{Example}[section]

\usepackage{caption}

\begin{document}

\maketitle
\begingroup
\makeatletter
\renewcommand\thefootnote{}
\renewcommand\@makefntext[1]{%
  \noindent #1%
}
\makeatother
\footnotetext{%
\textbf{Funding:} The second and third authors acknowledge research support by
the Polish National Science Centre, grant no.~2024/53/B/ST1/00703.
}
\endgroup

% REQUIRED
\begin{abstract}
We study discrete-time Markov Decision Processes (MDPs) on finite state-action spaces and analyze the stability of optimal policies and value functions in the long-run discounted risk-sensitive objective setting. Our analysis addresses robustness with respect to perturbations of the risk-aversion parameter and the discount factor, the emergence of ultimate stationarity, and the interaction between discounted and averaged formulations under suitable mixing assumptions. We further investigate limiting regimes associated with vanishing discount and vanishing risk sensitivity, and discuss the role of Blackwell-type stability properties in the discounted setting. Finally, we provide numerical illustrations that highlight the intrinsic non-stationarity of optimal discounted risk-sensitive policies.
\end{abstract}


 %REQUIRED

\section{Introduction}
Markov Decision Processes (MDPs) with discounted or averaged cumulative rewards over an infinite time horizon have attracted increasing attention, partly due to their central role in reinforcement learning; see \cite{Bor2010,Jai2022,FeiYanCheWan2021,PorCavCru2023,BisBor2023,GraPet2023} and references therein. While in classical MDP theory the key objective is to maximise the expected cumulative reward, say $Z$, the risk-sensitive extension replaces the expectation $\bE[Z]$ by the entropic utility 
\[
\mu^{\gamma}(Z):=\tfrac{1}{\gamma}\ln \bE\left[e^{\gamma Z} \right]
\]
with a pre-defined risk-aversion parameter $\gamma\neq 0$. 
For basic results on optimal risk-sensitive control in various settings, we refer to \cite{howard1972risk,Jaq1973,Whi1990,BraCavFer1998,DiMSte1999,CavFer2000,CavHer2009,chavez2015continuity,CavHer2017,CavCru2017,Ste2023,Ste2024,PitSte2024} and references therein; see also \cite{AsiJas2017,bauerle2014more,BisBor2023,bauerle2024markov} for extensions and alternatives.

The entropic modification of the objective criterion is motivated by risk–reward utility theory (see \cite{BiePli2003}) and leads to a risk-sensitive criterion whose Taylor expansion around $\gamma=0$ is given by $\mu^{\gamma}(Z) \approx \bE[Z] + \tfrac12 \gamma \textnormal{Var}(Z)$, highlighting its ability to capture variability-sensitive preferences; in the limit one recovers the classical risk-neutral setup. Among various risk-reward criteria, the entropic utility is particularly popular due to its strong time-consistency and tractable analytical properties, see \cite{KupSch2009,BiePli2003}.

Although discounted and averaged risk-sensitive problems are structurally different, they are linked through non-trivial limit properties, such as those related to Blackwell optimality or ultimate stationarity; while the existence theory for these types of problem is well established, the joint dependency of the optimal policy on the risk-sensitivity parameter $\gamma$ and the discounting parameter $\beta$ requires more in-depth studies.  In a companion paper \cite{BauPitSte2024}, we analyzed the stability of optimal policies only w.r.t.\ the risk-sensitivity parameter $\gamma$ in the {\em long-run stochastic control framework} with an averaged risk-sensitive criterion.  In this paper, we focus on the discounted framework and comment on the intricate links between averaged and discounted formulations and the interplay between the risk-sensitivity parameter $\gamma$ and the discount factor $\beta$. One could argue that the long-run averaged setup inherently focuses on time-stationary Markov policies, in contrast to the discounted setting discussed in this paper, where time-stationarity is in general no longer present. As in \cite{BauPitSte2024}, we restrict here to finite state and action spaces.

Beyond their intrinsic interest in stochastic control, the questions addressed in this paper are also closely related to recent developments in risk-sensitive reinforcement learning. Discounted formulations dominate reinforcement learning practice, while average-reward criteria often serve as benchmarks for long-run performance, see \cite{Put1994,SutBar2018}.
A rigorous understanding of the structural properties of optimal discounted policies, including their potential non-stationarity and their behavior under vanishing-discount limits, is essential for interpreting discounted approximations used in learning algorithms. In this context, the notion of ultimate stationarity provides a natural theoretical explanation for why, even in risk-sensitive settings where optimal discounted policies are generally non-stationary, stationary behavior may eventually emerge. Related structural phenomena have recently been exploited in the analysis of non-stationary and risk-sensitive reinforcement learning schemes, including settings where policy stabilization or eventual stationarity plays a key role, see  \cite{LeeDinLeeJinLavSoj2023,KumBhaSanKavHem2019,FeiYanWan2021,DinJinLav2023}.

Our main results can be summarised as follows. First, we study the optimal policy dynamics for small risk aversion and discuss the concept of \emph{ultimate stationarity}, originally considered in \cite{Jaq1972}. We independently show that for an arbitrary discount factor $\beta$ and the risk-sensitivity parameter $\gamma\neq 0$ there is an ultimately stationary optimal policy for the risk-sensitive problem, see Theorem~\ref{th:jaquette} and Theorem~\ref{th:jaquette2}. Although this result was already stated in \cite{Jaq1973}, the corresponding proof contains several steps that are difficult to verify in detail, which motivated us to provide an alternative and, hopefully, clearer argument. Building on this result, we also present a series of properties of optimal discounted policies. For instance, we show that for arbitrary $\beta$ there always exists a small $\gamma$ interval around zero on which an optimal stationary Markov policy can be found, i.e.\ it behaves like in the classical risk-neutral setup, and also study the behavior of the turnpike, i.e.\ the time at which the policy becomes ultimately stationary. Furthermore, we show that the small aversion $(\beta,\gamma)$-parameter subspace can be decomposed into a finite family of closed subsets that intersect at the boundary and include optimal ultimately stationary policies. We further show that the boundary of the region on which a single stationary policy stays optimal is lower semicontinuous in the discount factor, give conditions under which it is locally real analytic, and provide an example showing that it need not be continuous, see Corollary~\ref{cor:boundary}, Proposition~\ref{pr:gamma.continuity}, and Example~\ref{ex:4}. The decomposition of the $(\beta,\gamma)$-subspace is obtained under a mild non-degeneracy assumption on that boundary.

Second, under suitable mixing assumptions, we study the relation between discounted and averaged setups by analyzing the limiting behavior of optimal discounted policies and value functions as $\beta \to 1$. Here, we utilise the concept of Blackwell optimality and also recap the results of the basic vanishing discount framework, aligning the results of \cite{BauPitSte2024} with the discounted setup.

Third, we analyze the properties of the value functions and optimal policies in the vanishing risk-sensitivity regime, that is, when $\gamma \to 0$. We show that under mixing, the discounted Bellman operator is a span-contraction, and provide a discount factor uniform bound on the corresponding value function, reinforcing the result obtained in~\cite[Proposition 4.2]{BauPitSte2024}, see Theorem~\ref{th:contraction.discounted} and Proposition~\ref{pr:new.constant}. We also show that the distance between value functions in the risk-sensitive and risk-neutral setups can be controlled and show how risk-averse policies interact with risk-neutral policies for small aversion, see Proposition~\ref{pr:span-norm.fun} and Theorem~\ref{th:jaquette33}.  We also establish connections to moment optimal policies, illustrating why risk-sensitive analysis remains relevant even when the ultimate objective is to maximise the risk-neutral criterion.

The remainder of the paper is organised as follows. Section~\ref{S:introduction} collects preliminaries and comments on discounted Bellman policies. Section~\ref{S:Bellman} introduces the discounted Bellman equation and basic existence and regularity results for value functions and optimal policies.
Section~\ref{S:small.aversion} studies small risk-aversion dynamics, including ultimate stationarity, stationary optimal policies for small risk sensitivity, and turnpike-type notions. In Section~\ref{S:risk.average} we analyze the vanishing-discount limit and relate the behavior of optimal policies to Blackwell optimality.
Section~\ref{S:u-optimality} investigates the vanishing risk-sensitivity regime and its implications for mixing and moment optimality. Finally, Section~\ref{S:examples.all} presents numerical illustrations. Technical results related to the entropic utility are collected in Appendix~\ref{A:entropic}.

\section{Preliminaries}\label{S:introduction}

Let us consider a discrete-time MDP on a finite state-action space and follow the notation from~\cite{BauPitSte2024}, see also \cite{BauRie2011}. We use $E=\{x_1,\ldots,x_k\}$ and $U=\{\tilde u_1,\ldots,\tilde u_l\}$, where $k,l\in\bN$, to denote the state space and the action space, respectively. For simplicity, we assume that for any $x\in E$ the action set is full, that is, every action is always available; for $a\in U$, we use $P^{a}:=[p_a(x_i,x_j)]_{i,j=1}^{k}$ to denote the corresponding transition probability matrix.

We set a canonical measurable space $(\Omega,\cF)=(E^{\bN},(2^E)^{\bN})$ with filtration $(\cF_n)_{n\in\bN}$ where $\cF_n:=(2^E)^n\times\{\emptyset,E\}^\bN$, for $n\in\bN$, and use $\Pi$ to denote the space of all measurable policies on this space, that is, functions $\pi=(a_n)_{n\in\bN}$, where $a_n\colon\Omega\to U$ is $\cF_n$-measurable. Given $\pi\in\Pi$, we use $\bP^{\pi}$ to denote the corresponding controlled probability on $(\Omega,\cF)$ and use $\bP_x^{\pi}$ to denote the controlled probability for the canonical state process $(X_n)$ controlled by $\pi$ and starting at $x\in E$; the respective expectation operators are given by $\bE^{\pi}$ and $\bE^{\pi}_{x}$. For consistency, we also set $\bP^{a}(x,A):=\sum_{y\in A}p_a(x,y)$, $A\in 2^E$, for any $x\in E$ and action $a\in U$; note that $\bP^{a}(x,\cdot)$ is a one-step transition kernel.

We use $\Pi''$ to denote the set of all Markov policies and associate elements of $\Pi''$ with families of functions $(u_n)_{n\in\bN}$, $u_n\colon E\to U$, which induce policies $\pi=(a_n)_{n\in\bN}\in \Pi$, where $a_n=u_n(X_n)$. In addition, we use $\Pi'$ to denote the set of stationary Markov policies and link elements of $\Pi'$ with functions $u\colon E\to U$, which constitute policies $\pi=(a_n)_{n\in\bN}\in \Pi$, where $a_n=u(X_n)$. Note that $\Pi'\subset \Pi''\subset \Pi$ and $|\Pi'|=l^k<\infty$. 

For a fixed running reward function $c\colon E\times U\to\bR$, risk-sensitivity parameter $\gamma\in\bR$, discount factor $\beta\in (0,1)$, control strategy $\pi=(a_i)_{i\in\bN}\in\Pi$, and starting point $x\in E$ we consider {\it the $\beta$-discounted long run risk-sensitive criterion} that is given by
\begin{equation}\label{eq:RSC.disc}
J_\gamma(x,\pi;\beta):=
\begin{cases}
\tfrac{1}{\gamma}\ln\bE_x^{\pi}\left[e^{\gamma\sum_{i=0}^{\infty}\beta^i c(X_i,a_i)}\right], & \textrm{if }\gamma\neq 0,\\
\bE_x^{\pi}\left[\sum_{i=0}^{\infty}\beta^i c(X_i,a_i)\right], & \textrm{if } \gamma=0.
\end{cases}
\end{equation}
Here, $\gamma<0$ corresponds to the risk-averse regime, while $\gamma>0$ denotes the risk-seeking regime, see~\cite{Whi1990}.
For completeness, we summarise selected properties of the value function \eqref{eq:RSC.disc} in the finite setup, when the Markov policy is fixed; see also \cite{Jaq1976} and \cite{Jaq1972}. %analyticity of $\beta\to J_{\gamma}(x,\pi;\beta)$ for stationary policies is stated in page 48, line 6 in \cite{Jaq1976}.

\begin{proposition}[Properties of the policy value function]\label{pr:stability}
Fix $x\in E$ and $\pi\in\Pi$. Then
\begin{enumerate}
\item[a)] For $\beta\in (0,1)$, the function $\gamma\to J_{\gamma}(x,\pi;\beta)$ is bounded and non-decreasing.
\item[b)] The function $(\beta,\gamma)\to J_{\gamma}(x,\pi;\beta)$ is jointly real analytic on $(0,1)\times \bR$.
\end{enumerate}
\end{proposition}

\begin{proof}
Fix $x\in E$ and $\pi\in\Pi$. For any fixed $\beta\in (0,1)$, boundedness and monotonicity follows directly from the fact that the entropic utility measure is a monetary risk measure that is monotone with respect to the risk-aversion parameter $\gamma\in \bR$.  We get  $\tfrac{-\|c\|}{1-\beta} \leq J_{\gamma}(x,\pi;\beta)\leq \tfrac{\|c\|}{1-\beta}$; see Appendix~\ref{A:entropic}. Let $c_i := c(X_i,a_i)$, for $i\in\bN$. For brevity, we fix $S(\beta) := \sum_{i=0}^{\infty} \beta^i c_i$ and $M(\beta,\gamma):=\bE^{\pi}_{x}\left[e^{\gamma S(\beta)}\right]$. Noting the uniform bounds $|\gamma S(\beta)|\leq \frac{|\gamma|\cdot\|c\|}{1-\beta}$ and $|(S(\beta))^n|\leq \frac{\|c\|^n}{(1-\beta)^n}$, $n\in\bN$, and using the Taylor series expansion, we get
\begin{align*}
\textstyle M(\beta,\gamma)
&=  \textstyle \bE^{\pi}_{x}\left[\sum_{n=0}^\infty \frac{\gamma^n}{n!} (S(\beta))^n\right]=\sum_{n=0}^\infty \frac{\gamma^n}{n!} \bE^{\pi}_{x}\left[(S(\beta))^n\right]\\
&\textstyle =\sum_{n=0}^\infty \frac{\gamma^n}{n!}\bE^{\pi}_{x}\left[\sum_{k=0}^\infty a_{n,k}\beta^k\right] =\sum_{n,k\in\bN}\frac{\mathbb{E}_x^\pi[a_{n,k}]}{n!}\,\gamma^n\beta^k,
\end{align*}
for $\beta\in [0,1)$ and $\gamma\in\bR$, where $a_{n,k}:= \sum_{\substack{i_1,\ldots,i_n\in\bN_0\\ i_1+\ldots+i_n=k}} c_{i_1}\cdot \ldots\cdot c_{i_n}$ satisfies $\bE_x^{\pi}[|a_{n,k}|]<\binom{k+n-1}{n-1} \|c\|^n$, for $n,k\in\bN$; the last equality follows, e.g., from the dominated convergence theorem, as the series are absolutely convergent. The convergence of the series is uniform on any compact rectangle $R=[0,\beta_0)\times (-\gamma_0,\gamma_0]$, and we get that the mapping $M: (0,1)\times\bR\to\bR$ is jointly real-analytic. Now, since $M$ is strictly positive and the logarithm is real-analytic on $(0,\infty)$, the composition of $\log$ and $M$ is also real-analytic, and so is the mapping 
$(\beta,\gamma)\mapsto (\log M(\beta,\gamma))/\gamma$ on $(0,1)\times (\bR\setminus\{0\})$. Furthermore, the singularity at $\gamma=0$ is removable as entropic utility converges to expectation for $\gamma\to 0$; see Appendix~\ref{A:entropic}.
\end{proof}

\section{Discounted Bellman equation and optimal discounted policies}\label{S:Bellman}

Our main goal is to maximise the $\beta$-discounted long run risk-sensitive criterion given in \eqref{eq:RSC.disc}. We consider the optimal value function $w^{\beta}: E\times \bR \to \bR$ for problem \eqref{eq:RSC.disc} that is given by
\begin{equation}\label{eq:w.beta.32}
\textstyle w^{\beta}(x,\gamma):=\sup_{\pi\in\Pi}J_\gamma(x,\pi;\beta).
\end{equation}
This function is linked to the discounted risk-sensitive Bellman equation
\begin{equation}\label{eq:Bellman2}
\textstyle w^{\beta}(x, \gamma)=
\begin{cases}
\max_{a\in U}\left[c(x,a)+\frac{1}{\gamma}\ln\sum_{y\in E}e^{\gamma \beta w^{\beta}(y,\beta \gamma)}\bP^a(x,y)\right] & \textrm{if } \gamma\neq 0,\\
\max_{a\in U}\left[c(x,a)+\sum_{y\in E}\beta w^{\beta}(y,0)\bP^a(x,y)\right] & \textrm{if } \gamma= 0,
\end{cases}
\end{equation}
that can be reformulated using the discounted Bellman operator, that is, we can rewrite \eqref{eq:Bellman2} as 
\begin{equation}\label{eq:Bellman3}
\textstyle w^{\beta}(x, \gamma)=T^{\beta}\textstyle w^{\beta}(x, \gamma), 
\end{equation}
where
\begin{equation}\label{eq:Top}
\textstyle T^\beta h(x,\gamma):=
\begin{cases}
\max_{a\in U}\left[c(x,a) + \tfrac{1}{\gamma}\ln\sum_{y\in E} e^{\gamma\beta h(y,\beta\gamma)}\mathbb{P}^a(x,y)\right], & \textrm{if } \gamma\neq 0,\\
\max_{a\in U}\left[c(x,a)+\beta\sum_{y\in E} h(y,0)\bP^a(x,y)\right] ,& \textrm{if } \gamma= 0,
\end{cases}
\end{equation}
for  $h\in \bC(E\times \bR)$. Note that, for $\gamma\neq 0$, \eqref{eq:Bellman2} and \eqref{eq:Bellman3} could be seen as a system of equations defined on the risk-aversion grid $(\beta^i\gamma)_{i\in\bN}$ that is encoded as a second parameter of function $w^{\beta}$ in order to allow the recursive discounting scheme.\footnote{For simplicity, we consider the whole parameter space, rather than discretized parameter grid. Also, note that given the initial parameter $\gamma >0$ (resp. $\gamma<0$) the function $w^{\beta}$ in \eqref{eq:w.beta.32} can be effectively defined in $E\times (0,\gamma]$ (resp. $E\times [\gamma,0)$).}  Furthermore, for a fixed $\gamma\in \bR$, we consider the Markov policy $\hat\pi^{\beta}:=(\hat{u}^{\beta}_0,\hat{u}^{\beta}_1,\ldots)\in\Pi''$ induced by \eqref{eq:Bellman2} and given by
\begin{equation}\label{eq:u.non.optimal}
\textstyle \hat u^{\beta}_i(x):=
\begin{cases}
\argmax_{a\in U}[c(x,a) + \frac{1}{\gamma\beta^i}\ln ( \sum_{y\in E} e^{\gamma\beta^{i+1} w^\beta(y,\gamma \beta^{i+1})}\mathbb{P}^a(x,y))],& \textrm{if } \gamma\neq 0\\
\argmax_{a\in U}[c(x,a)+\sum_{y\in E}\beta w^{\beta}(y,0)\bP^a(x,y)] & \textrm{if } \gamma= 0,
\end{cases}
\end{equation}
for $i\in\bN_0$. If required, to ensure the uniqueness of the choice in \eqref{eq:u.non.optimal}, we always take the action with the smallest index from the subset of maximisers in $U=\{\tilde u_1,\ldots,\tilde u_l\}$. The link between \eqref{eq:w.beta.32}, \eqref{eq:Bellman2}, and \eqref{eq:u.non.optimal} is given in the following theorem; $\Pi^{*}(\beta,\gamma)$ is used to denote the set of optimal Markov policies for \eqref{eq:w.beta.32}  for all (starting points) $x\in E$.
\begin{theorem}\label{pr:exist2}
Fix $\beta\in (0,1)$. Then, the function $w^{\beta}$ defined in \eqref{eq:w.beta.32} solves the Bellman equation~\eqref{eq:Bellman2}. Furthermore, for any $\gamma\in \bR$, we have $\hat\pi^{\beta}\in\Pi^*(\beta,\gamma)$, that is, the policy 
$\hat\pi^{\beta}$ defined in \eqref{eq:u.non.optimal} is optimal for \eqref{eq:RSC.disc}.
\end{theorem}

\begin{proof}
Fix $\beta\in (0,1)$. For brevity, we show the proof only for $\gamma\neq 0$; the proof for $\gamma=0$ is similar. Define two sequences of functions $(\underline{w}^{\beta}_n)_{n\in\bN}$ and $(\overline{w}^{\beta}_n)_{n\in\bN}$, where  $\underline{w}^{\beta}_n\colon E\times \bR \to \bR$ and $\overline{w}^{\beta}_n\colon E\times \bR \to \bR$, for $n\in\bN$, are defined recursively by
\[
\begin{cases}
\underline{w}^{\beta}_0(x,\gamma):= -\frac{\|c\|}{1-\beta},\\
\underline{w}^{\beta}_{n+1}(x,\gamma):=T^\beta \underline{w}^{\beta}_n(x,\gamma),
\end{cases}
\qquad\textrm{and}\qquad
\begin{cases}
\overline{w}^{\beta}_0(x,\gamma):=\frac{\|c\|}{1-\beta},\\\overline{w}^{\beta}_{n+1}(x,\gamma):=T^\beta \overline{w}^{\beta}_n(x,\gamma).
\end{cases}
\]
It is easy to check that
\begin{align}
\underline{w}^{\beta}_n(x,\gamma) &=\sup_{\pi\in\Pi} \tfrac{1}{\gamma}\ln \mathbb{E}_x^{\pi}\left[e^{\sum_{i=0}^{n-1} \gamma \beta^i c(X_i,a_i)-\gamma\beta^n \tfrac{\|c\|}{1-\beta}}\right],\label{por1}\\
 \overline{w}^{\beta}_n(x,\gamma) &=\sup_{\pi\in\Pi} \tfrac{1}{\gamma}\ln \mathbb{E}_x^{\pi}\left[e^{\sum_{i=0}^{n-1} \gamma \beta^i c(X_i,a_i)+ \gamma\beta^n \tfrac{\|c\|}{1-\beta}}\right],\label{por2}
\end{align}
which implies that for any $(x,\gamma)\in E\times\bR$ and $n\in\bN$ we get
\begin{equation}\label{eq:3:approx1}
\tfrac{-\|c\|}{1-\beta}\leq \underline{w}^{\beta}_n(x,\gamma)\leq w^{\beta}(x,\gamma)\leq \overline{w}^{\beta}_n(x,\gamma)\leq \tfrac{\|c\|}{1-\beta},
\end{equation}
where $w^\beta$ is given in \eqref{eq:w.beta.32}. 
%Furthermore, for any $n\in\bN$, we have 
%\begin{align}
%\overline{w}^{\beta}_{n+1}(x,\gamma) &\leq \sup_{\pi\in\Pi} \tfrac{1}{\gamma}\ln \mathbb{E}_x^{\pi}\left[e^{\sum_{i=0}^{n-1} \gamma \beta^i c(X_i,a_i)+\gamma\beta^n\|c\|+ \gamma\beta^{n+1} \tfrac{\|c\|}{1-\beta}}\right]=  \overline{w}^{\beta}_n(x,\gamma),\label{eq:w.beta.new1}\\
%\underline{w}^{\beta}_{n+1}(x,\gamma) &\geq \sup_{\pi\in\Pi} \tfrac{1}{\gamma}\ln \mathbb{E}_x^{\pi}\left[e^{\sum_{i=0}^{n-1} \gamma \beta^i c(X_i,a_i)-\gamma\beta^n\|c\|- \gamma\beta^{n+1} \tfrac{\|c\|}{1-\beta}}\right]=  \underline{w}^{\beta}_n(x,\gamma).\label{eq:w.beta.new2}
%\end{align}
Combining \eqref{por1}, \eqref{por2}, and \eqref{eq:3:approx1}, and noting that, for $n\in\bN$, we have
\begin{equation}
0\leq \overline{w}^{\beta}_n(x,\gamma)-\underline{w}^{\beta}_n(x,\gamma)\leq \tfrac{2\beta^n\|c\|}{1-\beta},
\end{equation}
we get that the sequences $(\underline{w}^{\beta}_n)_{n\in\bN}$ and $(\overline{w}^{\beta}_n)_{n\in\bN}$ have the same limit equal to $w^\beta$. The optimality of $\hat\pi^{\beta}$  follows in a standard way using a verification procedure of Bellman equation \eqref{eq:Bellman2}; see \cite{BauRie2011} for details.
\end{proof}

For $\gamma\neq 0$, the risk-sensitive optimal Markov policy $\hat\pi^{\beta}$ defined in \eqref{eq:u.non.optimal} can be genuinely non-stationary, that is, for $i,j\in\bN$, such that $i\neq j$, we can get $\hat u^{\beta}_i\not\equiv \hat u^{\beta}_j$, see Section~\ref{S:examples} for illustration. On the other hand, for $\gamma=0$, the policy is stationary as for $i,j\in\bN$, we get $\hat u^{\beta}_i\equiv \hat u^{\beta}_j$ by construction. This shows a structural difference between the risk-sensitive and the risk-neutral discounted setup, as in the former case no optimal stationary policy may exist. For completeness, we also summarise some simple properties of the value function.

\begin{proposition}[General Properties of the optimal discounted policies and value functions]\label{th:stability2}
Fix $x\in E$. Then
\begin{enumerate}
\item[a)] For $\beta\in (0,1)$, the function $\gamma\to w^{\beta}(x,\gamma)$ is bounded and non-decreasing.
\item[b)] The function $(\beta,\gamma)\to w^{\beta}(x,\gamma)$ is continuous on $(0,1)\times \bR$.
\end{enumerate}
\end{proposition}

\begin{proof}
The boundedness and monotonicity follow immediately from the properties of the entropic utility measure; see Proposition~\ref{pr:entropic.monotone}. The proof of b) follows directly from the proof of Theorem~\ref{pr:exist2} and Proposition \ref{pr:stability} b).
%by taking into account that $T^\beta$ transforms the set of continuous functions into itself and noting that $w^\beta$ is the iteration limit for $T^{\beta}$.

\end{proof}

\section{Small risk-aversion dynamics and ultimate stationarity}\label{S:small.aversion} One expects that for small risk aversion, the risk-sensitive discounted policy will begin to resemble the risk-neutral policy and will admit some form of stationarity. To formalise this intuition, let us recall the concept of {\it ultimate stationarity} introduced in \cite{Jaq1973}, and investigate how it interacts with risk sensitivity.

We say that policy $\pi\in \Pi''$ is {\it ultimately  stationary}, if there exists $N\in\bN$ and $u_0,u_1,\ldots,u_{N-1},u\in\Pi'$ such that $\pi=(u_0,u_1,\ldots,u_{N-1},u,u,\ldots)$. In other words, the Markov policy $\pi$ is ultimately stationary if it becomes a stationary Markov policy after a fixed number of initial steps. The main goal of this section is to show that for any $\gamma\in\bR\setminus\{0\}$ and $\beta\in (0,1)$, there exists an ultimate stationary policy that is optimal for the discounted risk-sensitive problem. 

We first show that for any discount factor $\beta\in (0,1)$ and risk aversion $\gamma\in\bR\setminus\{0\}$ sufficiently close to zero one can find an optimal stationary policy. In this paper, we provide a modified and expanded version of the original proof from~\cite{Jaq1973}, where this result was initially formulated.  Note that parts of the proof presented here are structurally different from the ones in~\cite{Jaq1973}, as they are based on a direct analysis of Bellman's equation. We refer to Remark~\ref{rem:Jaquette.proof.problem}, where we provide a comment on some subtle intricacies and challenges that we encountered while analyzing the original proof in~\cite{Jaq1973}.

\begin{theorem}[Existence of an optimal stationary discounted policy for a small risk aversion]\label{th:jaquette}
Fix $\beta\in (0,1)$. Then, there exists $\gamma_\beta<0$ (resp. $\gamma_\beta>0$) and $u\in\Pi'$ such that $u\in \Pi^{*}(\beta,\gamma)$, for $\gamma\in [\gamma_\beta,0]$ (resp. $\gamma\in [0,\gamma_\beta]$).
\end{theorem}

%%%%%%%
%%%%%%%%%
%%%%%%%%%%

\begin{proof}

Fix $\beta<1$ and consider $\gamma<0$; the proof for $\gamma>0$ is analogous. For transparency, we split the proof into five steps.

\medskip 

\noindent {\it Step 1)} We show that there exists $\gamma_{0}<0$ such that the function 
\begin{equation}\label{eq:Hx}
H_x^{u,\tilde u}(\gamma):=J_{\gamma}(x,u,\beta)-J_{\gamma}(x,\tilde u,\beta), \quad\gamma\in\bR
\end{equation}
does not change the sign on $[\gamma_{0},0]$, for any $x\in E$ and $u,\tilde u\in\Pi'$. First, note that for any $u\in\Pi'$ and $x\in E$, the function $\gamma\to \bE^{u}_x[e^{\gamma\sum_{i=0}^{\infty}\beta^i c(X_i,u(X_{i}))}]$ is analytic on $\bR$ as the moment generating function of a bounded random variable, see \cite{Cur1942}. Consequently, for $x\in E$ and $u,\tilde u\in\Pi'$, there exists $\gamma_x^{u,\tilde u}<0$ such that the function
\begin{equation}\label{eq:hxuu}
h_x^{u,\tilde u}(\gamma)
:=\bE^{u}_x\left[e^{\gamma\sum_{i=0}^{\infty}\beta^i c(X_i,u(X_{i}))}\right]-\bE^{\tilde u}_x\left[e^{\gamma\sum_{i=0}^{\infty}\beta^i c(X_i,\tilde u(X_i))}\right]
\end{equation}
does not change the sign for $\gamma\in [\gamma_x^{u,\tilde u},0]$; the function $h_x^{u,\tilde u}$ is also analytic, so that it either has a finite number of zeros in any closed interval containing 0 or is equal to zero in this interval, see \cite[Corollary 1.2.7]{KraPar2002}. Second, fix
\[
\gamma_0:=\max_{x\in E, u,\tilde u\in\Pi'}\gamma_x^{u,\tilde u},
\]
and recall that sets $E$ and $\Pi'$ are finite. Consequently, we get that $\gamma_0<0$, and the function $h_x^{u,\tilde u}$ does not change the sign on $[\gamma_{0},0]$ for any $x\in E$ and $u,\tilde u\in\Pi'$. This directly implies that \eqref{eq:Hx} does not change the sign as well; the signs of $h_x^{u,\tilde u}$ and $H_x^{u,\tilde u}$ are opposite.

\medskip 

\noindent {\it Step 2)} We show that there exists a stationary policy $u^*\in\Pi'$ that is optimal among all stationary policies for $\gamma \in [\gamma_0,0]$. First, note that for fixed $x\in E$, the function given in \eqref{eq:Hx} could be used to constitute a partial order on the set of Markov policies. Namely, we get that $(\Pi',\preceq_x)$ is the partially ordered set, where $u \preceq_x \tilde u$ corresponds to $H^{u,\tilde u}_x(\gamma)\leq 0$, for $\gamma\in [\gamma_{0},0]$. As the set $\Pi'$ is finite, there exists the greatest element in $(\Pi',\preceq_x)$, say $u^x\in\Pi'$. Second, fix
\[
u^*(x):=u^x(x),\quad  x\in E,
\]
and note that $u^*\in\Pi'$. Then, for $x\in E$, and $u\in\Pi'$, we get
\begin{equation}\label{eq:best.u.star}
J_{\gamma}(x,u^*,\beta)\geq J_{\gamma}(x,u,\beta),\quad \textrm{for } \gamma \in [\gamma_{0},0],
\end{equation}
that is, $u^*$ is optimal within the set of Markov stationary policies for $\gamma \in [\gamma_{0},0]$; this follows directly from the fact that $H_x^{u^*, u}(\cdot)\geq 0$ on $[\gamma_{0},0]$. Indeed, using the standard recursive scheme, one can show that
\begin{align*}
J_\gamma(x,u, \beta) & \leq J_\gamma(x,u^x, \beta)\\
&\textstyle=\frac{1}{\gamma}\ln\sum_{y\in E}\exp\left({\gamma \left[c(x, u^x(x))+\beta J_{\beta\gamma}(y,u^x, \beta)\right]}\right)\bP^{u^x(x)}(x,y)\\
&\textstyle\leq \frac{1}{\gamma}\ln\sum_{y\in E}\exp\left({\gamma \left[c(x, u^x(x))+\beta J_{\beta\gamma}(y,u^y, \beta)\right]}\right)\bP^{u^x(x)}(x,y)\\
&\textstyle= \frac{1}{\gamma}\ln\sum_{y\in E}\sum_{z\in E}\exp\big(\gamma \big[c(x, u^x(x))+\beta c(y, u^y(y))\\
&\textstyle\phantom{\leq}\qquad\qquad\qquad\qquad+\beta^2 J_{\beta^2\gamma}(z,u^y, \beta)\big]\big)\bP^{u^y(y)}(y,z)\bP^{u^x(x)}(x,y)\\
&\textstyle\leq \frac{1}{\gamma}\ln\sum_{y\in E}\sum_{z\in E}\exp\big(\gamma \big[c(x, u^x(x))+\beta c(y, u^y(y))\\
&\textstyle\phantom{\leq}\qquad\qquad\qquad\qquad+\beta^2 J_{\beta^2\gamma}(z,u^z, \beta)\big]\big)\bP^{u^y(y)}(y,z)\bP^{u^x(x)}(x,y)\\
&\textstyle\leq \ldots\\
&\textstyle\leq J_\gamma(x,u^*, \beta).
\end{align*}

\medskip

\noindent {\it Step 3)} We show that there exists $\gamma_{1}<0$ such that $u^*$ is optimal for $\gamma\in [\gamma_{1},0]$, when confronted with non-stationary Markov policies in which we change the first action, that is, policies given by
\begin{equation}\label{eq:one.stat}
\textstyle \pi_{(u,1)}:=(u,u^*,u^*,u^*,\ldots),
\end{equation}
for $u\in\Pi'$. Following a reasoning similar to the one presented in Step 1), we know that for any $x\in E$ and $u\in\Pi'$ there exists $\gamma^u_x < 0$ such that the function
\[
H_x^{u}(\gamma):=J_{\gamma}(x,u^*,\beta)-J_{\gamma}(x,\pi_{(u,1)},\beta), \quad\gamma\in\bR,
\]
does not change the sign on $[\gamma^u_x,0]$. Let us fix
\[
\gamma_1:=\max_{x\in E, u\in \Pi'}\gamma_x^u \vee \gamma_0
\]
and show that $H_x^{u}(\gamma)\geq 0$ for $\gamma\in [\gamma_1,0]$, for any $x\in E$ and $u\in\Pi'$. On the contrary, taking into account that $H_{x}^{u}$ does not change the sign on $(\gamma_1,0)$, as $\gamma_0\leq \gamma_1$, assume that there exists $\bar x\in E$ and $\bar u\in\Pi'$ such that
\begin{equation}\label{eq:contr}
H_{\bar x}^{\bar u}(\gamma)< 0,\qquad\textrm{for } \gamma\in (\gamma_1,0),
\end{equation}
and let
\[
\tilde u(x):=
\begin{cases}
u^*(x) & \textrm{ if } x\neq \bar x,\\
\bar u(x) &\textrm { if } x=\bar x.
\end{cases}
\]
Consider a sequence of non-stationary policies given by
\begin{equation}\label{eq:pi.n}
\textstyle \pi_{(\tilde u,n)}:=(\, \underbrace{\tilde u,\tilde u,\ldots,\tilde u}_\text{n times}, u^*,u^*,u^*,\ldots),\quad n\in\bN.
\end{equation}
Then, using standard induction, we can show that, for any $n\in\bN$, we have
\begin{equation}\label{eq:ind}
J_\gamma(\cdot ,\pi_{(\tilde u,n)},\beta)\geq J_\gamma(\cdot ,\pi_{(\tilde u,n-1)},\beta) \qquad  \gamma\in (\gamma_1,0).
\end{equation}
For $n=1$, \eqref{eq:ind} follows directly from \eqref{eq:contr}. Now, assume that \eqref{eq:ind} holds for $n\in\bN$. Then, for any $x\in E$ and $\gamma\in (\gamma_1,0)$, noting that $\beta\gamma\in (\gamma_1,0)$, we have
\begin{align*}
J_\gamma(x ,\pi_{(\tilde u,n)},\beta)& = c(x,\tilde u(x))+\frac{1}{\gamma}\ln\sum_{y\in E}e^{\gamma\beta J_{\beta\gamma}(y,\pi_{(\tilde u,n-1)}, \beta)}\bP^{\tilde u(x)}(x,y)\\
& \leq c(x,\tilde u(x))+\frac{1}{\gamma}\ln\sum_{y\in E}e^{\gamma\beta J_{\beta\gamma}(y,\pi_{(\tilde u,n)}; \beta)}\bP^{\tilde u(x)}(x,y)\\
& = J_\gamma(x ,\pi_{(\tilde u,n+1)},\beta),
\end{align*}
which concludes the induction. Next, it is easy to show that
\begin{equation}\label{eq:limitt}
J_\gamma(\cdot ,\pi_{(\tilde u,n)},\beta) \to J_\gamma(\cdot ,\tilde u,\beta),\quad\textrm{ as } n\to\infty,
\end{equation}
using the standard continuity arguments; note that the total cumulative discounted reward received after $n$-steps can be bounded by $\pm \frac{\beta^n\|c\|}{1-\beta}$, while the first $n$ rewards for policies  $\pi_{(\tilde u,n)}$ and $\tilde u$ are the same. Combining \eqref{eq:ind} with \eqref{eq:limitt} we get
\[
J_\gamma(\bar x ,\tilde u,\beta)= \lim_{n\to\infty}J_\gamma(\bar x ,\pi_{(\tilde u,n)},\beta) \geq J_\gamma(\bar x ,\pi_{(\tilde u,1)},\beta)>J_\gamma(\bar x ,u^*,\beta),
\]
which contradicts \eqref{eq:best.u.star} as stationary policy $\tilde u$ outperforms $u^*$ on $(\gamma_1,0)$ and $\gamma_0<\gamma_1$.

\medskip

\noindent {\it Step 4)} We show that $u^*$ is optimal among all non-stationary Markov policies for $\gamma\in [\gamma_{1},0]$. Let us fix arbitrary policy $\tilde \pi=(\tilde u_1,\tilde u_2,\ldots)$ and the corresponding policies
\[
\textstyle \pi_{(\tilde \pi,n)}:=(\, \underbrace{\tilde u_1,\tilde u_2,\ldots,\tilde u_n}_\text{n times}, u^*,u^*,u^*,\ldots),\quad n\in\bN.
\]
Using Step 3), and applying similar reasoning as in the proof of Step 1), one can show that
\begin{equation}\label{eq:ind2}
J_\gamma(\cdot ,\pi_{(\tilde \pi,n-1)},\beta)\geq J_\gamma(\cdot ,\pi_{(\tilde \pi,n)},\beta) \qquad  \gamma\in (\gamma_1,0), \quad n\in\bN.
\end{equation}
In a nutshell, the first $n-1$ rewards for policies $\pi_{(\tilde \pi,n-1)}$ and $\pi_{(\tilde \pi,n)}$ are the same, and from Step 3) we know that it is better to apply $u^*$ than  $\tilde u_n$ in step $n$, if we later apply only $u^*$; note that for $n=1$, \eqref{eq:ind2} follows directly from Step 3). Also, following similar reasoning as in Step 3), we know that
\begin{equation}\label{eq:ind3}
J_\gamma(\cdot ,\pi_{(\tilde \pi,n)},\beta)\to J_\gamma(\cdot ,\tilde \pi,\beta)\quad \textrm{as } n\to\infty,
\end{equation}
Combining \eqref{eq:ind2} and \eqref{eq:ind3} and noting that $J_\gamma(\cdot ,\pi_{(\tilde \pi,0)},\beta)=J_\gamma(\cdot ,u^*,\beta)$, we conclude the proof of Step 4).

\medskip

\noindent {\it Step 5)} To conclude the proof, we need to show that $u^*$ is optimal among all policies for $\gamma\in [\gamma_{\beta},0]$, for $\gamma_\beta:=\gamma_1$. This is a direct consequence from the Bellman equation, which indicates that there exists a non-stationary Markov policy which maximises the discounted problem value, and Step 4).
\end{proof}

\begin{remark}[Comment on the original proof of Theorem~\ref{th:jaquette}]\label{rem:Jaquette.proof.problem}
As we already mentioned, Theorem~\ref{th:jaquette} was originally formulated and proved in \cite{Jaq1973}. When analysing the proof in \cite{Jaq1973}, we found many subtle challenges and informal reasoning which made the proof very hard to follow. This is why we have decided to present a full and modified version of the original proof in this paper; due to the best of our knowledge, the original proof was not recreated in any other paper or book. To better understand our motivation, let us provide an example. While the thesis in Lemma 3.3 from \cite{Jaq1973} is easy to show, the proof presented in \cite{Jaq1973} is unclear. In particular, we were unable to recreate two key inequalities posed in the proof while the proof mentioned that they follow directly from the definition. It should be also noted that the proof reasoning presented in Lemma 3.3 is later extrapolated to Lemma 3.4 without a proper comment. In fact, we were able to show that the reasoning presented in Lemma 3.4 is wrong, in a sense that inequality (3.9) is not true for $\lambda_0$ understood as maximum over $\lambda_0$s for all stationary policies and starting points as claimed in the paper, see definition of $\gamma_0$ in our proof for reference. We were able to remediate this by defining $\gamma_1$ in our version of the proof, but the reasoning presented in \cite{Jaq1973} does not mention this subtlety. 
\end{remark}

Naturally, it is of interest to choose $\gamma_\beta$ in Theorem~\ref{th:jaquette} as small as possible (or as large as possible when $\gamma_\beta>0$). Thus, let us define
 \begin{eqnarray}
        \gamma(\beta) &:=& \inf\{ \gamma_0 <0\colon \mbox{ there exists a decision rule } u \mbox{ such that } \nonumber\\
        &&  \qquad u \in \Pi' \mbox{ and } u\in\Pi^*(\beta,\gamma)  \mbox{ for all } \gamma \in [\gamma_0,0] \}.\label{eq:gamma.beta} 
    \end{eqnarray} 
For completeness, in analogy to \eqref{eq:gamma.beta}, we also set
\begin{eqnarray}
        \gamma_+(\beta) &:=& \sup\{ \gamma_0 >0\colon \mbox{ there exists a decision rule } u \mbox{ such that } \nonumber\\
        &&  \qquad u \in \Pi' \mbox{ and } u\in\Pi^*(\beta,\gamma)  \mbox{ for all } \gamma \in [0,\gamma_0] \},\label{eq:gamma.beta+} 
    \end{eqnarray} 
keeping in mind that our primary focus is on the risk-averse case. We can write $\gamma(\beta)$ in a different form. For $u\in\Pi', x\in E$ and $a\in U$, let $\pi^{x,a,u}:= (v^{x,a},u,u,\ldots)$ where $v^{x,a}:E\to U$ is any decision rule satisfying $v^{x,a}(x)=a.$ Its values at states different from $x$ are irrelevant when the initial state is $x.$ Further define the one-step-deviation difference
\[
D_{u,x,a}(\beta,\gamma) := J_\gamma(x,u,\beta)-J_\gamma(x,\pi^{x,a,u},\beta), \quad (\beta,\gamma)\in (0,1)\times \mathbb{R}.
\]
For $u\in \Pi'$ and $\beta\in (0,1)$ define
\begin{equation}\label{def:rhou}
\rho_u(\beta):= \inf\{r: D_{u,x,a}(\beta,\gamma) \ge 0 \mbox{ for all } x\in E, a\in U, \gamma \in[r,0]\}
\end{equation} 
where the infimum of the empty set is understood to be $+\infty$, so that $\rho_u(\beta)\in[-\infty,0]\cup\{+\infty\}$, and we use the convention $[-\infty,0]:=(-\infty,0]$. We finally obtain that
\[
\gamma(\beta)= \min_{u\in\Pi'} \rho_u(\beta).
\]
The last equation holds because a stationary policy $u$ is optimal for all $\gamma\in[r,0]$ if and only if $D_{u,x,a}(\beta,\gamma)\geq0$ for all $x\in E$, $a\in U$, and $\gamma\in[r,0]$, i.e.\ we only have to compare a stationary policy against policies which differ in the first decision rule. This follows by the same argument as in Steps 4)--5) of the proof of Theorem~\ref{th:jaquette}, since $[r,0]$ is invariant under $\gamma\mapsto\beta\gamma$ and hence the one-step comparison propagates along the risk-aversion grid $(\beta^{i}\gamma)_{i\in\bN}$.
%    
%The last equation is true since we have shown in the proof of Theorem~\ref{th:jaquette} that a stationary policy $u$ is optimal for every $\gamma \in [r,0]$ if and only if $D_{u,x,a}(\beta,\gamma) \geq 0.$, i.e. we only have to compare a stationary policy against policies which differ in the first decision rule.
%
Using this definition, we get the following result.

\begin{corollary}\label{cor:boundary}
    Let $\gamma(\beta) $ be given by \eqref{eq:gamma.beta}.
    \begin{enumerate}
    \item[a)] For every $\beta\in (0,1)$ it holds that $\gamma(\beta)\in[-\infty,0)$ and there exists a $u_\beta \in \Pi'$ such that $u_\beta \in \Pi^*(\beta,\gamma)$ for every $\gamma \in [\gamma(\beta),0].$ In particular, if $\gamma(\beta)>-\infty$, the maximal interval is closed at its left endpoint.
    \item[b)] The mapping $\beta\mapsto \gamma(\beta)$ is lower semicontinuous.
    \end{enumerate}
    The analogous statements hold for the risk-seeking interval to the right of zero.
\end{corollary}

\begin{proof}
    Part a). The fact that $\gamma(\beta)\in [-\infty,0)$  follows directly from Theorem~\ref{th:jaquette} and the existence of  $u_\beta \in \Pi'$ such that $u_\beta \in \Pi^*(\beta,\gamma)$ for every $\gamma \in (\gamma(\beta),0]$ from the definition of $\gamma(\beta).$ What remains to show in part a) is the optimality at the left endpoint $\gamma(\beta).$ 

    Now assume that $(r_n)$ with $ r_n<0$ is such that $r_n\downarrow \gamma(\beta)$ where  $\gamma(\beta)<r_n$ (with $r_n\downarrow-\infty$ if $\gamma(\beta)=-\infty$). Then $D_{u_\beta,x,a}(\beta,\gamma)\ge 0$ for all $x,a$ and every $\gamma \in [r_n,0].$ The functions $(\beta,\gamma) \mapsto D_{u,x,a}(\beta,\gamma)$ are continuous, in fact jointly real analytic, because they are differences of fixed policy value functions. Letting $n\to\infty$ we obtain $D_{u_\beta,x,a}(\beta,\gamma)\ge 0$ for all $\gamma \in [\gamma(\beta),0].$ Thus, we also obtain optimality at the endpoint of the interval. This proves part a).

    Part b). Fix $u\in\Pi'$ and $r\leq 0$. Define 
    $$ m_{u,r}(\beta) := \min_{x\in E, a\in U, \gamma \in [r,0]} D_{u,x,a}(\beta,\gamma).$$
    Because $E,U$ are finite, the interval $[r,0]$ is compact and the functions $D_{u,x,a}$ are jointly continuous, the mapping $\beta \mapsto m_{u,r}(\beta)$ is continuous. Moreover,
    $$ \{ \beta\in(0,1) : \rho_u(\beta) \le r\} = \{ \beta\in(0,1) : m_{u,r}(\beta)\ge 0\}.$$
    The set on the r.h.s.\ is closed. Thus, every lower-level set of $\rho_u$ is closed, and hence $\rho_u$ is lower semicontinuous. Since $\gamma(\beta)$ is the finite minimum of  functions $\rho_u$ the statement is shown.
\end{proof}

The preceding semicontinuity statement can be strengthened at regular boundary points. 

\begin{proposition}[Continuity of $\beta\to \gamma(\beta)$]\label{pr:gamma.continuity}
Fix $\beta_0\in(0,1)$ and $\gamma_0:=\gamma(\beta_0)$.
Assume that there exist a stationary policy $u_0\in\Pi'$, a state $x_0\in E$, and an action $a_0\in U$, such that 
\begin{enumerate}
    \item the policy $u_0$ is the unique maximal radius policy, that is, we have $\rho_{u_0}(\beta_0) <\min_{u\neq u_0}\rho_u(\beta_0)$;
    \item the point $(x_0,a_0)$ is the unique active constraint at the left endpoint, that is, we have $a_0\neq u_0(x_0)$, $D_{u_0,x_0,a_0}(\beta_0,\gamma_0)=0$, $D_{u_0,x_0,a_0}(\beta_0,\gamma)>0$, for $\gamma\in (\gamma_0,0]$, and $D_{u_0,x,a}(\beta_0,\gamma)>0$, for $\gamma\in [\gamma_0,0]$ and $(x,a)$ such that $a\neq u_0(x)$ and $(x,a)\neq (x_0,a_0)$.

    \item the point $(x_0,a_0)$ is a simple active zero, that is, $\partial_\gamma D_{u_0,x_0,a_0}(\beta_0,\gamma_0) \neq0$.
\end{enumerate}
Then there exist a neighborhood $V$ of $\beta_0$ and a real-analytic function $\widehat\gamma : V\to(-\infty,0)$ such that $\widehat\gamma(\beta_0)=\gamma_0$ and $\widehat\gamma(\beta)=\gamma(\beta)$, for $\beta\in V$.
In particular, $\gamma(\beta)$ is continuous and real analytic in a neighborhood of $\beta_0$.
\end{proposition}

\begin{proof}
    Set $F(\beta,\gamma):=D_{u_0,x_0,a_0}(\beta,\gamma).$
The function $F$ is jointly real analytic. From the assumptions, we know that
$F(\beta_0,\gamma_0)=0$, and $\partial_\gamma F(\beta_0,\gamma_0)\neq0.$ The real-analytic implicit-function theorem (see, e.g., \cite[Section 2]{KraPar2002}) yields neighborhoods
$V$ of $\beta_0$ and $W$ of $\gamma_0$, together with a unique real-analytic function
$\widehat\gamma : V\to W$ such that
$F(\beta,\widehat\gamma(\beta))=0$
and $\widehat\gamma(\beta_0)=\gamma_0.$

Since $F(\beta_0,\gamma)>0$ for
$\gamma\in(\gamma_0,0]$, and the zero at $\gamma_0$ is simple, $F$ changes sign at $\gamma_0$. Shrinking $V$ and $W$, if necessary, we obtain
$$ F(\beta,\gamma)>0
\quad
\text{for }
\widehat\gamma(\beta)<\gamma\leq0,
\mbox{ and }
F(\beta,\gamma)<0
\quad
\text{for }
\gamma<\widehat\gamma(\beta),$$
whenever $\gamma$ is sufficiently close to $\widehat\gamma(\beta).$
For every other $(x,a)$, joint continuity of the deviation functions and our assumptions imply $D_{u_0,x,a}(\beta,\gamma)>0$ for all $\beta \in V, \gamma \in [\hat\gamma(\beta),0].$ Thus, $\rho_{u_0}(\beta) = \hat\gamma(\beta)$ for $\beta \in V.$

Finally, because $u_0$ is the unique maximiser at the stability radius at $\beta_0$, there exists $\varepsilon>0$ such that $\rho_u(\beta_0) \ge \rho_{u_0}(\beta_0)+3\varepsilon$  for every $u\neq u_0.$ By the lower semicontinuity of $\rho_u$ (see previous proof), after shrinking $V$ we have $\rho_u(\beta) \ge \rho_{u_0}(\beta_0)+2\varepsilon$  for every $u\neq u_0$ and $ \beta\in V.$ On the other hand $\rho_{u_0}(\beta) =\hat\gamma(\beta) \le \rho_{u_0}(\beta_0)+\varepsilon$  for $\beta$ sufficiently close to $\beta_0.$ Thus, $u_0$ remains the unique policy attaining the maximal radius and $\gamma(\beta)=\hat\gamma(\beta).$
\end{proof}

In the next result, we show how to use Theorem~\ref{th:jaquette} to show existence of optimal ultimately stationary Markov policy for the generic discounted problem. We also refer to \cite[Theorem 1]{Jaq1976} that is based on the results obtained in \cite{Jaq1973}.

\begin{theorem}[Existence of ultimately stationary policy]\label{th:jaquette2}
For any $\gamma\in\bR$ and $\beta\in (0,1)$ there exists optimal Markov policy $\pi\in \Pi^{*}(\beta,\gamma)$ that is ultimately stationary.
\end{theorem}

\begin{proof}
The proof follows directly from the combination of Theorem~\ref{pr:exist2},  Theorem~\ref{th:jaquette} and the recursive discounted problem formulation. For $\gamma=0$ the proof is immediate as $\hat\pi^{\beta}$ defined in \eqref{eq:u.non.optimal} is stationary. Fix $\gamma\in\bR\setminus\{0\}$ and $\beta\in (0,1)$. Noting that $\beta^i\gamma \to 0$, as $i\to\infty$, and using Theorem~\ref{th:jaquette}, we get that there exists $i_0\in\bN$ and $u\in\Pi'$ such that $\bar \pi:=(u,u,\ldots,)\in\Pi^{*}(\beta,\beta^i\gamma)$ for $i\ge i_0$. From Theorem~\ref{pr:exist2}, we also know that $\hat\pi^{\beta}=(\hat{u}^{\beta}_0,\hat{u}^{\beta}_1,\ldots)\in \Pi^{*}(\beta,\gamma)$, for $\hat\pi^{\beta}$ is defined in \eqref{eq:u.non.optimal}. Let us define ultimately stationary policy $\tilde\pi\in\Pi''$ given by
\begin{equation}\label{eq:ult.stat.proof1}
\textstyle \tilde \pi=(\tilde u_0,\tilde u_1,\ldots):=(\, \underbrace{\hat{u}^{\beta}_0,\ldots,\hat{u}^{\beta}_{{i_0}-1}}_\text{$i_{0}$ times}, u,u,u,\ldots).
\end{equation}
To show that policy \eqref{eq:ult.stat.proof1} is optimal for \eqref{eq:RSC.disc} it is sufficient to note that it solves the sequence of Bellman equations \eqref{eq:Bellman2} on sensitivity parameter grid ($\beta^i\gamma)_{i=0}^{\infty}$, that is, for any $i\in\bN$, we get
\begin{equation}\label{eq:iteration.proof1}
\textstyle w^{\beta}(x, \beta^i\gamma)=c(x,\tilde u_i(x))+\frac{1}{\beta^i\gamma}\ln\sum_{y\in E}e^{\gamma \beta^{i+1} w^{\beta}(y,\beta^{i+1} \gamma)}\bP^{\tilde u_i(x)}(x,y).
\end{equation}
Indeed, for $i=0,1,\ldots,i_0-1$, \eqref{eq:iteration.proof1} follows directly from the definition of policy $\hat\pi^{\beta}$, see \eqref{eq:u.non.optimal}, while for $i\ge i_0$, using direct recursion, for $x\in E$, we get
\begin{align*}
w^{\beta}(x, \beta^i\gamma) & =\sup_{\pi\in\Pi}J_{\beta^i\gamma}(x,\pi;\beta)=J_{\beta^i\gamma}(x,\bar\pi;\beta)= \tfrac{1}{\beta^i\gamma}\ln\bE_x^{\bar\pi}\left[e^{\gamma\beta^i\sum_{k=0}^{\infty}\beta^k c(X_k,u(X_k))}\right]\\
& = c(x,u(x))+\tfrac{1}{\beta^i\gamma}\ln\bE_x^{\bar\pi}\left[e^{\gamma\beta^i\sum_{k=1}^{\infty}\beta^k c(X_k,u(X_k))}\right]\\
& =c(x,u(x))+\tfrac{1}{\beta^i\gamma}\ln\sum_{y\in E}\bE_y^{\bar\pi}\left[e^{\gamma\beta^{i+1}\sum_{k=0}^{\infty}\beta^k c(X_k,u(X_k))}\right]\bP^{u(x)}(x,y)\\
& = c(x,u(x))+\tfrac{1}{\beta^i\gamma}\ln\sum_{y\in E}e^{\gamma\beta^{i+1} w^{\beta}(y,\beta^{i+1}\gamma)}\bP^{u(x)}(x,y),
\end{align*}
which concludes the proof, as $\tilde u_i(x)=u(x)$, for $i\ge i_0$.
\end{proof}

\begin{remark}[Ultimate stationarity on infinite state spaces]
The equivalent of Theorem~\ref{th:jaquette2} is no longer true for denumerable state space, that is, one can show that there might be no ultimately stationary policy for the discounted problem if the state space is non-finite. For details, we refer to Example 4 in \cite{BraCavFer1998}.
\end{remark}

For any $\beta\in (0,1)$ and $\gamma\in\bR$, we use $N(\beta,\gamma)\in\bN$ to denote the {\it turnpike}, that is, an earliest ultimate stationarity threshold. In other words, $N(\beta,\gamma)$ denotes the smallest iteration after which there exists an optimal policy which is stationary from the $N(\beta,\gamma)$-th stage onwards. Let us now analyse the behaviour of the function $(\beta,\gamma)\to N(\beta,\gamma)$. 

\begin{corollary}\label{cor:N}
    Fix $\beta\in (0,1)$. Then, $\gamma \to N(\beta,\gamma)$ 
    is locally piecewise constant, with only finitely many jumps on compact intervals.
 %   is piecewise constant.
    %, increasing in $(-\infty,0]$ and decreasing in $[0,\infty)$.
\end{corollary}

\begin{proof}

%\NB{Obviously for all $\gamma \in [\gamma_0/\beta^n,\gamma_0/\beta^{n+1} )$ the turnpike $N(\beta,\gamma)$ is the same by definition. Monotonicity is also clear by construction.} \MP{NB, can you please elaborate on those two observations? I proved that $N$ is indeed piecewise constant, but the proof of monotonicity is not obvious to me -- see comment below. Also, I do not understand why turnpike is constant on the quoted intervals? }

Fix $\beta\in (0,1)$. For brevity, we focus on $[0,\infty)$; the proof for $(-\infty,0]$ is similar. %First, let us show that the mapping is piecewise constant. 
Using Theorem \ref{th:jaquette}, we know that there exists $\gamma_+(\beta)>0$, such that $N(\beta,\gamma)=1$ for $\gamma\in [0,\gamma_+(\beta)]$. From the iteration procedure embedded into the discounted Bellman equation we also get $N(\beta,\beta^n\gamma)=\max\{N(\beta,\beta^{n-1}\gamma)-1,1\}$
for any $n\in\bN$ and $\gamma\geq 0$. Those two observations imply that
$$
N(\beta,\gamma)\leq n\quad  \textrm{for } \gamma\in [0,\beta^{-n}\gamma_+(\beta)].
$$
Consequently, noting that the set of ultimately stationary Markov policies with turnpike smaller or equal to $n$ is finite and using the real-analytic property of the corresponding value functions in a similar way as in the proof of Theorem~\ref{th:jaquette} (see also Proposition~\ref{pr:stability}) we get that the function $\gamma\to N(\beta,\gamma)$ changes its value finitely many times on the interval $[0,\beta^{-n}\gamma_+(\beta)]$, for any $n\in\bN$. Thus, noting that $\beta^{-n}\gamma_+(\beta)\to\infty$ as $n\to\infty$, we get that $\gamma\to N(\beta,\gamma)$ changes its value finitely many times on every compact subinterval.
\end{proof}
%\MP{I am not sure about monotonicity. I started to draft the proof but it is not clear to me if this property is even true.}

%\MP{\noindent {\bf [DRAFT TEXT]:} Let $\gamma^{1}_\beta:=\inf\{\gamma\in\bR_{+}\colon N(\beta,\gamma)>1\}$. From Theorem \ref{th:jaquette} we know that $\gamma_{\beta}^1>0$. If $\gamma^{1}_\beta=\infty$, the proof is concluded, as $N(\beta,\cdot)$ is constant on $[0,\infty)$. Assume that $\gamma^{1}_\beta<\infty$. From Proposition~\ref{pr:stability}, i.e., the continuity of the value function, we know that $N(\beta,\gamma^1_{\beta})=1$.
%From the iteration procedure we get that $N(\beta,\gamma)\leq 2$ for $\gamma\in (\gamma_{\beta}^1,\tfrac{\gamma_{\beta}^1}{\beta}]$ as well as $N(\beta,\gamma)=2$ for $\gamma\in (\gamma_{\beta}^1,\gamma_{\beta}^2]$; last equality is due to real-analyticity of the value function ans the fact that there are finitely many ultimately stationary Markov policy with turnpike smaller of equal to two. }

%\noindent \MP{I am not sure how to proceed from that. In particular how are we sure that we dont get $N(\beta,\gamma)=1$ for some $\gamma$s in $(\gamma_{\beta}^2,\tfrac{\gamma_{\beta}^1}{\beta}]$. NB, can you elaborate?}
%\NB{Is it possible to say something about $N=N(\beta,\gamma)$ from which time onwards $u$ is optimal? This reminds me of Turnpike numbers in \cite{LewAna2019}. I think the proof implies the following: For a compact set $[\gamma_0,\gamma_1]\times [\beta_0,\beta_1]$ there are only a finite number of turnpike values $N=N(\beta,\gamma)$.}

We conclude this section with a comment on the stability of optimal policies.

\begin{proposition}[Finite partition of $(\beta,\gamma)$-plane]\label{pr:finite.partition}
    Fix $\varepsilon \in (0,\tfrac12)$, $\gamma_0>0$, and set $\Gamma := [\varepsilon,1-\varepsilon]\times [-\gamma_0,\gamma_0].$ Assume that $\underline{\rho} := \inf_{\beta \in [\varepsilon,1-\varepsilon] } \min\{-\gamma(\beta),\gamma_+(\beta) \}>0.$ Then there exists a finite collection $\mathcal{P}_0\subset \Pi''$ of ultimate stationary Markov policies such that for every $(\beta,\gamma)\in \Gamma$ at least one policy in $\mathcal{P}_0$ is optimal for the discounted problem. More precisely, for every $\pi\in \mathcal{P}_0$ define
    \begin{align}
        K_\pi := \{(\beta,\gamma) \in \Gamma : J_\gamma(x,\pi,\beta)\ge J_\gamma(x,\tilde\pi,\beta) \mbox{ for every } x\in E, \tilde\pi\in\mathcal{P}_0\}
    \end{align}
    Then:
    \begin{itemize}
        \item[a)] every $K_\pi$ is a closed subset of $\Gamma$;
        \item[b)] the sets $K_\pi$ cover $\Gamma$;
        \item[c)] If $\pi,\sigma\in\mathcal{P}_0$ are not value-equivalent, then $K_\pi\cap K_\sigma \subset \partial_\Gamma K_\pi \cup \partial_\Gamma K_\sigma$, where $\partial_\Gamma$ denotes the relative boundary in $\Gamma$; we call $\pi,\sigma\in\mathcal{P}_0$ value-equivalent if $J_\gamma(x,\pi,\beta)=J_\gamma(x,\sigma,\beta)$ for $x\in E$ and $(\beta,\gamma)\in\Gamma$. 
    \end{itemize}

\end{proposition}

\begin{proof}
   Part a). Set $\bar\beta := 1-\varepsilon$ and define $N_0 := \max\{0, \lceil\frac{\log(\gamma_0/\underline{\rho})}{-\log\bar\beta}\rceil\}.$ By construction $\bar\beta^{N_0}\gamma_0 \le \underline{\rho}.$
    Let $\mathcal{P}_0$ denote the set of all Markov policies of the form \linebreak  $\pi=(u_0,u_1,\ldots,u_{N_0-1},u,u,\ldots)$, where $u_0,u_1,\ldots,u_{N_0-1},u : E \to U$ are decision rules. Since state and action space are finite, $\mathcal{P}_0$ is finite. By Theorem \ref{th:jaquette2} and the construction of $\Gamma$ it is clear that for every $(\beta,\gamma)\in\Gamma$ there exists an optimal policy in $\mathcal{P}_0$.  

    Now fix $\pi\in \mathcal{P}_0$. For every $x\in E, \tilde\pi\in \mathcal{P}_0$ the function 
    \begin{align*}
        (\beta,\gamma)\mapsto J_\gamma(x,\pi,\beta)-J_\gamma(x,\tilde \pi,\beta)
    \end{align*}
    is jointly real-analytic on $(0,1)\times \mathbb{R}.$ Hence $K_\pi$ is a finite intersection of sets of the form
    \begin{align*}
        \{(\beta,\gamma) \in \Gamma : J_\gamma(x,\pi,\beta)-J_\gamma(x,\tilde\pi,\beta)\ge 0 \}.
    \end{align*}
    Thus, it follows that $K_\pi$ is closed.

    Part b). Fix $(\beta,\gamma) \in \Gamma.$ Theorem \ref{th:jaquette2} implies that there exist a $\pi^*\in \mathcal{P}_0$ such that $(\beta,\gamma) \in K_{\pi^*}.$

    Part c). Suppose that there exists a point $(\beta,\gamma)\in K_\pi \cap K_\sigma$ which belongs neither to $\partial_\Gamma K_\pi$ nor to $\partial_\Gamma K_\sigma.$ Thus, $(\beta,\gamma)$ belongs to the relative interior of both $K_\pi$ and $K_\sigma.$ Hence there exists a relatively open neighborhood $V$ of $(\beta,\gamma)$ such that $V\subset K_\pi\cap K_\sigma.$ Since $\varepsilon\in(0,\tfrac12)$ and $\gamma_0>0$, the relatively open set $V$ contains a non-empty open subset of $(0,1)\times\bR$, on which $J_\gamma(x,\pi,\beta)= J_\gamma(x,\sigma,\beta)$ for every $x\in E$. As these functions are jointly real analytic on the connected set $(0,1)\times\bR$, see Proposition~\ref{pr:stability} b), the identity theorem implies that they coincide on $(0,1)\times\bR$, contradicting the assumption that $\pi$ and $\sigma$ are not value-equivalent. This proves the statement.
\end{proof}
It should be noted that the assumption $\underline\rho>0$ is not automatically satisfied, see Section~\ref{sec:Exgamma} for an example.

\section{Vanishing discount and Blackwell optimality}\label{S:risk.average}

In this section, we examine the limiting behavior of optimal discounted policies and value functions as $\beta \to 1$, and recall the standard vanishing-discount framework. For brevity, we decided to provide only an overview of the key concepts and results, since they are more directly linked to the risk-sensitive average stochastic control problems studied in our companion paper, see~\cite{BauPitSte2024} for details and proofs.

In order to establish a link between the discounted setup and the standard risk-average setup, in which the optimal value does not depend on the starting point,  one needs to impose additional mixing properties on the underlying MDP. To this end, we consider the assumptions from \cite[Section 2.1]{BauPitSte2024}:

\begin{enumerate}
\item[(\namedlabel{A.1}{A.1})] (One-step mixing.) We have that
\[
\overline\Delta:=\sup_{x,x'\in E}\sup_{A\in 2^E}\sup_{a,a'\in U} \left|\bP^a(x,A)-\bP^{a'}(x',A)\right| <1.
\]

\item[(\namedlabel{A.2}{A.2})] (Multi-step transition equivalence.) There exists $N\in\bN$ such that
\[
K:=\sup_{x,x'\in E}\sup_{y\in E}\sup_{\pi\in\Pi^{''}} \frac{\bP^{(\pi,N)}(x,y)}{\bP^{(\pi,N)}(x',y)} <\infty,
\]
where $\bP^{(\pi,N)}$ denotes the $N$th iteration of the controlled transition kernel and conventions $\tfrac{0}{0}:=1$ as well as $\tfrac{1}{0}:=\infty$ are used.
\end{enumerate}
\medskip
For simplicity, we decided to assume a strong one-step mix which could be linked to the {\it uni-chain} condition, i.e.\ the existence of a unique (controlled) recurrent class. Combination of assumptions \eqref{A.1} and \eqref{A.2} is sometimes called {\it strong mixing} under which the span (semi-)norm of the value function for the discounted problem can be controlled; for $f\in C(E)$ we define $\|f\|_{\textnormal{sp}}:=\sup_{x,y\in E}[f(x)-f(y)]$.

\begin{proposition}[Span-norm value function bounds induced by the mixing conditions]\label{pr:w.disc.bound}
Assume \eqref{A.1} and \eqref{A.2}. Then, for any $\gamma\neq 0$, we have
\[
\sup_{\beta\in (0,1)}\|w^{\beta}(\cdot,\gamma)\|_{\textnormal{sp}}\leq N\|c\|_{\textnormal{sp}}+\tfrac{1}{|\gamma|}\ln K,
\]
where $N\in\bN$ and $K<\infty$ are the constants stated in \eqref{A.2}. 
\end{proposition}

For the proof, we refer to \cite[Proposition 4.2]{BauPitSte2024}. Notably the induced span-norm  bound is independent of $\beta$. Consequently, in the next step, one wants to let $\beta\to 1$ and study the limit behaviour of the value function. Of course, the problem needs to be properly normalised as otherwise this value will diverge. This is typically done by considering a unit-per-time reward increment. Namely, for any $n\in\bN$, let us fix a starting point $\bar z\in E$, and define the entropy increment function
\begin{equation}\label{eq:new.lamba.n}
\lambda_n^{\beta}(\gamma):= \beta^n w^{\beta}(\bar z,\gamma\beta^n)-\beta^{n+1}w^{\beta}(\bar z,\gamma\beta^{n+1}),
\end{equation}
which essentially measures the entropic difference between the last $n$ and last $n+1$ elements of the discounted sum under the risk aversion re-normalised to the initial parameter $\gamma$. Under mixing, when $\beta\to 1$, one would expect that for any $n\to\infty$ the increment would converge to a single value that can be in turn linked to the risk sensitive averaged reward per unit step; see \cite{CavFer2000} for more background. The corresponding value function could be simply defined as a difference between discounted criteria for different starting points, e.g., it can be defined as
\begin{equation}\label{eq:new.w.n}
\bar w_n^{\beta}(x,\gamma):=\beta^nw^{\beta}(x,\gamma\beta^n)-\beta^n w^{\beta}(\bar z,\gamma\beta^n).
\end{equation}
Let us now formulate the vanishing discount convergence result which quantifies the idea presented beforehand.

\begin{theorem}[Vanishing discount theorem]\label{th:vanishing.discount}
Assume \eqref{A.1} and \eqref{A.2}. Then, for any $\gamma\neq 0$, there exists a function $ w(\cdot,\gamma)\colon E\to \bR$ and a constant $\lambda(\gamma)\in\bR$, such that, for $n\in\bN$ and $x\in E$, we have
\[
\lim_{\beta\uparrow1} \bar w_n^{\beta}(x,\gamma)= w(x,\gamma)\quad \textrm{and}\quad \lim_{\beta\uparrow 1}\lambda_{n}^{\beta}(\gamma)=\lambda(\gamma),
\]
where $\lambda_{n}^{\beta}(\gamma)\in\bR$ and $\bar w_n^{\beta}(\cdot,\gamma)\colon E\to\bR$ are given in \eqref{eq:new.lamba.n} and \eqref{eq:new.w.n}.
\end{theorem}
For the proof, we refer to \cite[Theorem 4.4]{BauPitSte2024}. It should be noted that the limits introduced in Theorem~\ref{th:vanishing.discount} are solving the Bellman equation for the averaged problem, and could be used to recover the optimal value as well as optimal strategy for the stochastic control problem with long-run risk-sensitive averaged objective criterion
\begin{equation}\label{eq:RSC.averaged}
\tilde J_\gamma(x,\pi):=
\begin{cases}
\liminf_{n\to\infty}\tfrac{1}{n}\frac{1}{\gamma }\ln\bE_x^{\pi}\left[e^{\gamma\sum_{i=0}^{n-1}c(X_i,a_i)}\right],& \textrm{if } \gamma\neq 0,\\
\liminf_{n\to\infty}\frac{1}{n}\bE_x^{\pi}\left[\sum_{i=0}^{n-1}c(X_i,a_i)\right],& \textrm{if }\gamma=0.
\end{cases}
\end{equation}
For details, we again refer to our companion paper \cite{BauPitSte2024}.

The averaged problem is solved by a stationary Markov strategy. Let us now comment how it can be recovered from the discounted problem; note that we want to recover the strategy itself rather than the optimal value function, as in Theorem~\ref{th:vanishing.discount}. The existence of Blackwell optimal strategies had initially been studied under the risk-neutral setup, in which the following had been shown.

\begin{theorem}[Blackwell property in the risk-neutral setup]\label{th:blackwell}
Assume \eqref{A.1} and fix $\gamma=0$.  Then, there exists a stationary Markov policy $\pi\in\Pi'$ which satisfies the Blackwell property, i.e.\ there is $\beta_0\in (0,1)$ such that for any $\beta\in (\beta_0,1)$ we have $\pi\in \Pi^*(\beta,0)$. Furthermore the strategy is optimal for problem~\eqref{eq:RSC.averaged} with $\gamma=0$.
\end{theorem}

We refer to \cite{Bla1962} for details; note that assumption \eqref{A.1} has been added to Theorem~\ref{th:blackwell} for consistency with other results. A quick analysis of Theorem~\ref{th:blackwell} shows us that one cannot recover the same result for the risk-sensitive case. This is due to the structural difference between risk-neutral and risk-sensitive setups. For instance, while stationary policy always exists for the discounted risk-neutral problem, it might  not be the case for the risk-sensitive discounted problem. That said, one can still recover the optimal policy from the discounted setup.

\begin{theorem}[Blackwell property in the risk-sensitive setup]\label{th:blackwell.sensitive}
Assume \eqref{A.1} and \eqref{A.2}. Then, for any $\gamma\neq 0$ there is $\beta_0\in (0,1)$ such that for any $\beta\in(\beta_0,1)$ the stationary Markov policy $\hat u^{\beta}_0\in \Pi'$ defined in \eqref{eq:u.non.optimal} is optimal for the risk-sensitive averaged problem. 
\end{theorem}

For the proof, see \cite[Theorem 4.5]{BauPitSte2024}. Note that the stationary policy $\hat u^{\beta}_0$ in Theorem~\ref{th:blackwell.sensitive} is simply a first component of the discounted policy from the discounted setup. We have so far introduced two special types of stationary policies: the ultimate stationarity policy and the Blackwell optimal policy. It should be emphasised that those policies could be different as we later illustrate in 
Section~\ref{S:examples}.

\section{Vanishing risk-sensitivity and moment optimality}\label{S:u-optimality}
While in the previous case we considered the limit policies and value functions for $\beta \to 1$, in this section we focus on the limit policies and value functions for $\gamma\to 0$.  From Section~\ref{S:small.aversion}, we already know that for $\gamma$ sufficiently close to zero, there exists an optimal stationary
policy for the risk-sensitive discounted problem which is also optimal for the risk-neutral discounted problem, see Theorem~\ref{th:jaquette}. Now, let us focus on the interactions between the optimal policies and value functions for the risk-sensitive and the risk-neutral discounted criterion assuming mixing conditions, that is, \eqref{A.1} and \eqref{A.2}.

First, note that for the risk-neutral setup ($\gamma=0$), the discounted Bellman operator $T^{\beta}$ is a global contraction under one-step mixing.

\begin{theorem}\label{th:discounted.contraction.risk.neutral}
Assume \eqref{A.1} and fix $\gamma=0$. Then, for any $\beta\in (0,1)$, the discounted Bellman operator $T^{\beta}$ from \eqref{eq:Top} restricted to the set $\bC(E\times \{0\})$ is a span-norm contraction with respect to the first variable, that is, for any $f_1,f_2\in \bC(E\times \{0\})$ we get
\[
\|T^{\beta}f_1(\cdot,0)-T^{\beta}f_2(\cdot,0)\|_{\textnormal{sp}}\leq \overline\Delta \|f_1(\cdot,0)-f_2(\cdot,0)\|_{\textnormal{sp}},
\]
where $\overline\Delta<1$ is a shrinkage constant from   \eqref{A.1}.
\end{theorem}
We refer to \cite[Note 2 in Section 5.5]{HerLas1996} for the proof of Theorem~\ref{th:discounted.contraction.risk.neutral} and more details on span-contraction framework in the risk-neutral setup. A similar property is true for risk-sensitivity sufficiently close to zero. To formulate the result, we need to introduce a restricted family of functions 
\[
\bC_N(\gamma_0):=\{f\in \bC(E\times [\gamma_0,0))\colon f\equiv (T^{\beta})^N g \textrm{ for some } g\in \bC(E\times [\gamma_0,0))\},
\]
defined for any $\gamma_0<0$ and with $N\in\bN$ from \eqref{A.2}. In simple words, $\bC_N(\gamma_0)$ is the family of functions that is obtained after applying $N$ iterations of operator $T^{\beta}$ to a generic function from $\bC(E\times [\gamma_0,0))$; note that risk-sensitivity is also constrained to the set $[\gamma_0,0)$. This is done to ensure a proper level of mixing, so that transition equivalence based arguments could be utilised.

\begin{theorem}\label{th:contraction.discounted}
Assume \eqref{A.1} and \eqref{A.2}. Then, for any $\beta\in (0,1)$ and $\gamma_0<0$ the discounted Bellman operator $T^{\beta}$ restricted to the set $\bC_N(\gamma_0)$ is a span-norm contraction with respect to the first variable, that is, for any $\gamma\in [\gamma_0,0)$ and $f_1,f_2\in \bC_N(\gamma_0)$ we get
\[
\|T^{\beta}f_1(\cdot,\gamma)-T^{\beta}f_2(\cdot,\gamma)\|_{\textnormal{sp}}\leq \Delta \|f_1(\cdot,\gamma)-f_2(\cdot,\gamma)\|_{\textnormal{sp}}
\]
for a fixed $\Delta<1$. Furthermore, the constant $\Delta$ is independent of the choice of $\beta$.
\end{theorem}

\begin{proof}
For any $\beta\in (0,1)$, $f\in \bC(E\times \bR_{-})$, $\gamma_0<0$, $\gamma\in [\gamma_0,0)$, $d>0$, and $A\in\cP(E)$ we set a supplementary notation
\begin{align*}
\|f\|_{\gamma} & :=\sup_{\gamma'\in [\gamma,0)}\|\gamma'f(\cdot,\gamma')\|_{\textnormal{sp}},\nonumber\\
\Lambda(d) &:= \sup_{\beta \in (0,1)} \sup_{\gamma\in [\gamma_0,0)}\max_{x,x'\in E}\max_{a,a'\in U}\max_{A\in \cP(E)}\sup_{\|g_1\|_{\gamma}\leq d}\sup_{\|g_2\|_{\gamma}\leq d} \mathbb{H}^{x,x',a,a'}_{\gamma,\beta,g_1,g_2}(A),\nonumber\\
\mathbb{H}^{x,x',a,a'}_{\gamma,\beta,g_1,g_2}(A) & :=\nu_{x,a,\gamma,\beta,g_1}(A)-\nu_{x',a',\gamma,\beta,g_2}(A),\\\
\nu_{x,a,\gamma,\beta,f}(A) & :=\frac{\sum_{z\in A}e^{\gamma\beta f(z,\gamma\beta)}\bP^a(x,z)}{\sum_{y\in E}e^{\gamma\beta f(y,\gamma\beta)}\bP^a(x,y)}.
\end{align*}
For transparency, we split the remaining part of the proof into four steps.

\medskip

\noindent {\it Step 1)} Following min-max approach from the span-contraction framework  we get
\begin{equation}\label{eq:bound.good1}
\|T^{\beta}f_1(\cdot,\gamma)-T^{\beta}f_2(\cdot,\gamma)\|_{\textnormal{sp}}\leq \beta\|f_1(\cdot,\gamma)-f_2(\cdot,\gamma)\|_{\textnormal{sp}} \Lambda(\|f_1\|_{\gamma} \vee \|f_2\|_{\gamma}),
\end{equation}
for any functions $f_1,f_2\in \bC(E\times \bR_{-})$. For details, we refer, e.g., to \cite[Lemma 1]{PitSte2016}.
\medskip

\noindent {\it Step 2)} We show that for any fixed $d>0$ we get
\begin{equation}\label{eq:bound.good2}
\Lambda(d)<1.
\end{equation}
On the contrary, assume that it does not hold i.e.\ there is $\beta_n\to \beta$,
$\gamma_n\to \gamma$, $a_n\to a$, $x_n\to x$, $y_n\to y$, $A_n\in \cP(E)$ and functions $g_{1,n}$, $g_{2,n}\in\bC(E\times \bR_{-})$ such that $\|g_{1,n}\|\leq d$ and $\|g_{2,n}\|\leq d$ and 
$\nu_{x_n,a_n,\gamma_n,\beta_n,g_{1,n}}(A_n)\to 1$, while
$\nu_{x_n',a_n',\gamma_n,\beta_n,g_{2,n}}(A_n)\to 0$.
Since
\[
e^{-\|f\|_{\gamma}}P_x^a(A)\leq \nu_{x,a,\gamma,\beta,f}(A)\leq e^{\|f\|_{\gamma}}P_x^a(A)    
\]
% DIRECT PROOF OF THIS INEQUALITY:
%Noting that for any $\gamma_0<0$, $x,z\in E$, $a\in U$, $\gamma\in [\gamma_0,0)$, $\beta\in (0,1)$ and $f\in \bC_N(\gamma_0)$ we have
%\begin{align*}
%\nu_{x,a,\gamma,\beta,f}(z) & \leq \frac{\sup_{y\in E} e^{\gamma\beta f(z,\gamma\beta)}\bP^a(y,z)}{\inf_{y\in E}e^{\gamma\beta f(y,\gamma\beta)}}\leq e^{-\gamma\beta \|f\|_{\gamma_0}}\bP^a(x,y)\\
%\nu_{x,a,\gamma,\beta,f}(z) & \geq \frac{\inf_{y\in E} e^{\gamma\beta f(z,\gamma\beta)}\bP^a(y,z)}{\sup_{y\in E}e^{\gamma\beta f(y,\gamma\beta)}}\geq e^{\gamma\beta \|f\|_{\gamma_0}}\bP^a(x,y)
%\end{align*}
we have that $P_{x_n'}^{a_n'}(A_n)\to 0$ and $P_{x_n}^{a_n}(A_n^c)\to 0$. Therefore $P_{x_n'}^{a_n'}(A_n^c)\to 1$, and we get a contradiction with \eqref{A.1}.

\medskip

\noindent {\it Step 3)} We show that for any fixed $\gamma_0<0$ we get
\begin{equation}\label{eq:bound.good3}
\sup_{f\in \bC_N(\gamma_0)}\|f\|_{\gamma}<\infty.
\end{equation}
Fix $\gamma_0<0$ and $f\in \bC_N(\gamma_0)$. Recalling that $f\equiv (T^{\beta})^N g$ for some $g\in \bC(E\times [\gamma_0,0))$ and using \eqref{A.2}, for any $\gamma\in [\gamma_0,0)$ we get
\begin{align*}
\|\gamma (T^{\beta})^N g(\cdot,\gamma)\|_{\textnormal{sp}} & \leq \sup_{x,y\in E}\sup_{\pi\in\Pi}\ln \frac{\mathbb{E}_x^{\pi}\left[e^{\sum_{i=0}^{N-1} \gamma \beta^i c(X_i,a_i)+\gamma\beta^N g(X_n,\beta^N\gamma)}\right]}{\mathbb{E}_y^{\pi}\left[e^{\sum_{i=0}^{N-1} \gamma \beta^i c(X_i,a_i)+\gamma\beta^N g(X_n,\beta^N\gamma)}\right]}\nonumber\\
& \leq |\gamma| N \|c\|_{\textnormal{sp}} +\sup_{x,y\in E}\sup_{\pi\in\Pi}\ln \frac{\mathbb{E}_x^{\pi}\left[e^{\gamma\beta^N g(X_n,\gamma)}\right]}{\mathbb{E}_y^{\pi}\left[e^{\gamma\beta^N g(X_n,\gamma)}\right]}\nonumber\\
& \leq |\gamma_0| N \|c\|_{\textnormal{sp}} +\ln K,
\end{align*}
which concludes this step of the proof.

\medskip 
\noindent {\it Step 4)} Finally, combining \eqref{eq:bound.good1}, \eqref{eq:bound.good2}, and \eqref{eq:bound.good3} and setting $\Delta:=\Lambda(|\gamma_0| N \|c\|_{\textnormal{sp}} +\ln K)$, we get
\begin{equation}\label{eq:contr.3}
\|T^{\beta}f_1(\cdot,\gamma)-T^{\beta}f_2(\cdot,\gamma)\|_{\textnormal{sp}}\leq \Delta \|f_1(\cdot,\gamma)-f_2(\cdot,\gamma)\|_{\textnormal{sp}}.
\end{equation}
for any $f_1,f_2\in \bC_N(\gamma_0)$, which concludes the proof.
\end{proof}

\begin{remark}[Interaction between multi-step transition equivalence and small risk aversion dynamics]
In order to maintain conciseness and consistency with the frameworks in Section~\ref{S:risk.average} and Section~\ref{S:u-optimality}, we assume \eqref{A.1} and \eqref{A.2} in Theorem~\ref{th:contraction.discounted} and the subsequent results. However, it should be noted that assumption \eqref{A.2} can be dropped in cases where $\gamma$ is presumed to be sufficiently close to zero. For an in-depth examination of the relationship between \eqref{A.2} and the dynamics associated with small values of $\gamma$, we refer to \cite[Section 4]{PitSte2024}. Additionally, refer to \cite{PitSte2016} for an idea of alternative proof methodology for Step 3) of Theorem~\ref{th:contraction.discounted}.
\end{remark}

Using Theorem~\ref{th:contraction.discounted} we can show that the discounted optimal value is uniformly bounded in the span norm; this allows us to get rid of dependency on risk-aversion parameter when compared to the bound stated in Proposition~\ref{pr:w.disc.bound}.

\begin{proposition}[Span-norm value function uniform bound induced by the mixing conditions]\label{pr:new.constant}
Assume \eqref{A.1} and \eqref{A.2}, and fix $\gamma_0<0$. Then, for any $\gamma<0$, we have
\[
\sup_{\beta\in (0,1)}\|w^{\beta}(\cdot,\gamma)\|_{\textnormal{sp}}\leq \tfrac{3N+1}{1-\Delta}\|c\|_{\textnormal{sp}}+\tfrac{1}{|\gamma_0|}\ln K,
\]
where $N\in\bN$ and $K<\infty$ are the constants stated in \eqref{A.2}, and $\Delta\in (0,1)$ is the shrinkage constant from Theorem~\ref{th:contraction.discounted} corresponding to $\gamma_0$.
\end{proposition}

\begin{proof}

For $\gamma<\gamma_0$, the proof follows directly from Proposition~\ref{pr:w.disc.bound}, as $N\|c\|_{\textnormal{sp}}+\tfrac{1}{|\gamma|}\ln K\leq \tfrac{3N+1}{1-\Delta}\|c\|_{\textnormal{sp}}+\tfrac{1}{|\gamma_0|}\ln K$. Let us fix $\gamma\in [\gamma_0,0)$,  $\beta\in (0,1)$, and set $K_i(\cdot,\gamma):=(T^{\beta}0)^{N+i}(\cdot,\gamma)-(T^{\beta}0)^{N+i-1}(\cdot,\gamma)$, for $i\in\bN$. Then, for any $m\in\bN$, we get
\begin{align*}
\|w^{\beta}(\cdot,\gamma)\|_{\textnormal{sp}} &\leq \|w^{\beta}(\cdot,\gamma)-(T^{\beta}0)^{N+m}(\cdot,\gamma)\|_{\textnormal{sp}} +\|(T^{\beta}0)^{N+m}(\cdot,\gamma)\|_{\textnormal{sp}}\\
& \leq \tfrac{\beta^{N+m}}{1-\beta}\|c\|_{\textnormal{sp}}+\|(T^{\beta}0)^N(\cdot,\gamma) +\sum_{i=1}^{m}K_i(\cdot,\gamma)\|_{\textnormal{sp}}\\
&\leq \tfrac{\beta^m}{1-\beta}\|c\|_{\textnormal{sp}}+\|(T^{\beta}0)^N(\cdot,\gamma)\|_{\textnormal{sp}}+\sum_{i=1}^{m}\|K_i(\cdot,\gamma)\|_{\textnormal{sp}}\\
&\leq \tfrac{\beta^m}{1-\beta}\|c\|_{\textnormal{sp}}+N\|c\|_{\textnormal{sp}}+\sum_{i=1}^{\infty}\Delta^{i-1}\|K_1(\cdot,\gamma)\|_{\textnormal{sp}}\\
&\leq \tfrac{\beta^m}{1-\beta}\|c\|_{\textnormal{sp}}+N\|c\|_{\textnormal{sp}}+\frac{(2N+1)\|c\|_{\textnormal{sp}}}{1-\Delta}.
\end{align*}
Letting $m\to\infty$, we get $\|w^{\beta}(\cdot,\gamma)\|_{\textnormal{sp}} \leq \frac{(3N+1)\|c\|_{\textnormal{sp}}}{1-\Delta}$, for $\gamma\in [\gamma_0,0)$ and $\beta\in (0,1)$, which concludes the proof.
\end{proof}

In the next proposition we show that, under mixing, the distance between the discounted risk-sensitive and risk-neutral value functions can be controlled in the span-norm independently of the discount factor, for $\gamma$ sufficiently close to zero.

\begin{proposition}\label{pr:span-norm.fun}
Assume \eqref{A.1} and \eqref{A.2}. Then, we have
\[
\sup_{\beta\in (0,1)}\|w^{\beta}(\cdot,\gamma)-w^{\beta}(\cdot,0)\|_{\textnormal{sp}}\to 0,\quad \textrm{as } \gamma\to 0.
\]
\end{proposition}

\begin{proof}
Fix $\beta\in (0,1)$ and assume that $\gamma\in [\gamma_0,0)$, where $\gamma_0<0$ is from Theorem~\ref{th:contraction.discounted}. From the triangle inequality we get that for any $m\in\bN$ and $N\in\bN$ set in \eqref{A.2}, we have
\begin{align}
\|w^{\beta}(\cdot,\gamma)-w^{\beta}(\cdot,0)\|_{\textrm{sp}} & \leq \phantom{+}\|w^{\beta}(\cdot,\gamma)-(T^{\beta}0)^{m+N}(\cdot,\gamma)\|_{\textrm{sp}}\nonumber\\
&\phantom{\leq}+\|(T^{\beta}0)^{m+N}(\cdot,\gamma)-(T^{\beta}0)^{m+N}(\cdot,0)\|_{\textrm{sp}}\nonumber\\
& \phantom{\leq}+\|(T^{\beta}0)^{m+N}(\cdot,0)-w^{\beta}(\cdot,0)\|_{\textrm{sp}}\label{eq:newone.1}
\end{align}
Let us now provide upper bounds for each constituents of \eqref{eq:newone.1}. First, using Theorem~\ref{th:contraction.discounted} and Proposition~\ref{pr:new.constant}, for $\gamma\in [\gamma_0,0)$ and $m\in\bN$, we have
\begin{align}
\|w^{\beta}(\cdot,\gamma)-(T^{\beta}0)^{m+N}(\cdot,\gamma)\|_{\textrm{sp}} &= \|(T^{\beta}w^{\beta})^{m+N}(\cdot,\gamma)-(T^{\beta}0)^{m+N}(\cdot,\gamma)\|_{\textrm{sp}}\nonumber\\
& \leq \Delta^m \|(w^{\beta})^N(\cdot,\gamma)-(T^{\beta}0)^{N}(\cdot,\gamma)\|_{\textnormal{sp}}\nonumber\\
& = \Delta^m \|w^{\beta}(\cdot,\gamma)-(T^{\beta}0)^{N}(\cdot,\gamma)\|_{\textnormal{sp}}\nonumber\\
& \leq \Delta^m \left(\|w^{\beta}(\cdot,\gamma)\|_{\textnormal{sp}}+\|(T^{\beta}0)^{N}(\cdot,\gamma)\|_{\textnormal{sp}}\right)\nonumber\nonumber\\
&\leq \Delta^m\left(\left(\tfrac{3N+1}{1-\Delta}\|c\|_{\textnormal{sp}}+\tfrac{1}{|\gamma_0|}\ln K\right)+N\|c\|_{\textnormal{sp}}\right)\nonumber\\
&\leq \Delta^m(\tfrac{4N+1}{1-\Delta}\|c\|_{\textnormal{sp}}+\tfrac{1}{|\gamma_0|}\ln K)\nonumber\\
& =:a_1(m).
\label{eq:sp.vanish.1}
\end{align}
Second, using Hoeffding's lemma (see \eqref{eq:hoef} in  Appendix~\ref{A:entropic}), for $\gamma\in [\gamma_0,0)$ and $m\in \bN$,
\begin{align}
\|(T^{\beta}0)^{m+N}(\cdot,\gamma)-(T^{\beta}0)^{m+N}(\cdot,0)\|_{\textrm{sp}} & \leq \textstyle \frac{|\gamma|}{2}\left(\sum_{i=1}^{m+N}\|c\|_{\textnormal{sp}}\right)^2\nonumber\\
& \leq \tfrac{|\gamma|}{2}(m+N)^2\|c\|_{\textnormal{sp}}^2\nonumber\\
& =:a_2(m,\gamma).\label{eq:sp.vanish.2}
\end{align}
Third, following similar logic as in \eqref{eq:sp.vanish.1}, recalling that $w^{\beta}$ satisfies \eqref{eq:Bellman2}, and using \eqref{A.1}, for $\gamma\in [\gamma_0,0)$ and $m\in\bN$, we get
\begin{align}
\|(T^{\beta}0)^{m+N}(\cdot,0)-w^{\beta}(\cdot,0)\|_{\textrm{sp}} & \leq (\overline\Delta)^m(\|w^{\beta}(\cdot,0)\|_{\textrm{sp}}+N\|c\|_{\textrm{sp}})\nonumber\\
& \leq (\overline\Delta)^m\left(\frac{2\|c\|_{\textrm{sp}}}{1-\beta\overline\Delta}+N\|c\|_{\textrm{sp}}\right)\nonumber\\
& \leq (\overline\Delta)^m\left(\frac{2\|c\|_{\textrm{sp}}}{1-\overline\Delta}+N\|c\|_{\textrm{sp}}\right)\nonumber\\
& =:a_3(m).\label{eq:sp.vanish.3}
\end{align}
Combining \eqref{eq:newone.1}, \eqref{eq:sp.vanish.1}, \eqref{eq:sp.vanish.2}, and \eqref{eq:sp.vanish.3}, for $\gamma\in [\gamma_0,0)$ and $m\in\bN$, we get
\begin{equation}\label{eq:final.44}
\sup_{\beta\in (0,1)}\|w^{\beta}(\cdot,\gamma)-w^{\beta}(\cdot,0)\|_{\textrm{sp}}\leq a_1(m)+a_2(m,\gamma)+a_3(m).
\end{equation}
Letting $\gamma\to 0$ and noting that $a_2(m,\gamma)\to 0$, for any $m\in\bN$, we get
\begin{equation}\label{eq:final.4444}
\textstyle 0\leq \limsup_{\gamma\to 0}\left[\sup_{\beta\in (0,1)}\|w^{\beta}(\cdot,\gamma)-w^{\beta}(\cdot,0)\|_{\textrm{sp}}\right]\leq a_1(m)+a_3(m).
\end{equation}
Now, letting $m\to\infty$ and noting that $a_1(m)+a_3(m)\to 0$, we conclude the proof.
\end{proof}

Using Proposition~\ref{pr:span-norm.fun}, we can easily show that the limit statement in Theorem~\ref{th:jaquette} can be rephrased in the spirit of Theorem~\ref{th:blackwell.sensitive}, that is, the first constituent of any optimal discounted policy is also optimal for the risk-neutral discounted problem, for sufficiently small $\gamma$.

\begin{theorem}[Small risk aversion and risk-neutral optimality]\label{th:jaquette33}
Assume \eqref{A.1} and \eqref{A.2}. Then, for any $\beta\in (0,1)$ there is $\gamma_0<0$, such that the first constituent of the Markov policy defined in \eqref{eq:u.non.optimal}, for any $\gamma\in [\gamma_0,0]$, is a stationary Markov optimal policy for the risk-neutral discounted problem, that is, we have $\hat u^{\beta}_0\in \Pi^{*}(\beta,0)$.
\end{theorem}
\begin{proof}

While the proof of Theorem~\ref{th:jaquette33} could be deduced by modifying the proof of Theorem~\ref{th:jaquette}, we present an alternative proof based more directly on mixing. Fix $\beta\in (0,1)$ and assume that the statement of Theorem~\ref{th:jaquette33} does not hold. Then, there exists a decreasing sequence $(\gamma_k)_{k\in\bN}$, such that $\gamma_k\to 0$, as  $k\to\infty$, and a series of stationary Markov policies $(u_k)_{k\in\bN}$, $u_k\in\Pi'$, such that $u_k\not\in \Pi^*(\beta,0)$ and
\begin{equation}\label{eq:optimal.new.1}
u_k(x)\in \argmax_{a\in U}\left[c(x,a) + \frac{1}{\gamma_k}\ln \left( \sum_{y\in E} e^{\gamma_k \beta\, w^\beta(y,\gamma_k \beta)}\mathbb{P}^a(x,y)\right)\right],\quad x\in E.
\end{equation}
Then, since $E$ and $U$ are finite, there exists a stationary Markov policy $u\in \Pi'$, such that $u_{k}\equiv u$ for $k$ coming from an increasing subsequence of $\bN$, say $(k_m)_{m\in\bN}$. From Hoeffding's lemma (see Appendix~\ref{A:entropic}), fixing $\tilde \gamma_m :=\gamma_{k_m}$, for $m\in\bN$, we get
\[
\left|\frac{1}{\tilde \gamma_m}\ln \sum_{y\in E} e^{\tilde \gamma_m\beta w^\beta(y, \beta\tilde \gamma_m)}\mathbb{P}^a(x,y) -\sum_{y\in E}\beta w^\beta(y, \beta\tilde \gamma_m)\mathbb{P}^a(x,y) \right|\leq |\tilde \gamma_m|\frac{\|w^{\beta}(\cdot, \beta\tilde\gamma_m)\|_{\textnormal{sp}}^2}{8}.
\]
Consequently, using \eqref{eq:optimal.new.1}, for $x\in E$, we get
\[
w^{\beta}(x,\tilde\gamma_m)\leq c(x,u(x)) + \sum_{y\in E}\beta w^\beta(y, \beta\tilde \gamma_m)\mathbb{P}^{u(x)}(x,y)+|\tilde \gamma_m|\frac{\|w^{\beta}(\cdot, \beta\tilde\gamma_m)\|_{\textnormal{sp}}^2}{8}.
\]
Thus, letting $m\to\infty$ and using Proposition~\ref{pr:span-norm.fun}, we get
\[
w^{\beta}(x,0)\leq c(x,u(x)) +\sum_{y\in E}\beta w^\beta(y,0)\mathbb{P}^{u(x)}(x,y),\quad x\in E,
\]
which implies $u \in \Pi^{*}(\beta,0)$. This leads to the contradiction with the preimposed condition $u_k\not\in \Pi^*(\beta,0)$, $k\in\bN$, so that the statement of Theorem~\ref{th:jaquette33} must hold.
\end{proof}

Let us now focus on the optimal policies. As already noted, from Theorem~\ref{th:jaquette}, we get that for any $\beta\in (0,1)$ there exists a stationary Markov policy $\pi\in \Pi'$ such that $\pi\in\Pi^*(\beta,\gamma)$ for $\gamma\in [\gamma_{\beta},0]$. In particular, recalling that this policy is uniquely determined (up to degenerative symmetries, etc.), this indicates that one can recover discounted risk-neutral optimal policy by solving discounted risk-sensitive problem restricted to Markov stationary policies. This raises the question, if this is effectively worth doing, as the risk neutral problem is, in principle, easier to solve due to the linearity of the expectation operator. To answer this question, we need to recall the concept of {\it moment optimality} that has been initially introduced in \cite{Jaq1973}.

We say that a policy $\pi\in \Pi$ is {\it (risk-averse) moment optimal} for the $\beta$-discounted reward if $\pi$ is the maximal element for the lexicographic order induced by the moments vector, in which every odd moment is multiplied by (-1). 
More formally, we require that for any $\pi'\in \Pi$ and $x\in E$, there exists $n\in\bN$ such that
\begin{equation}\label{eq:moment.optimal}
\begin{cases}
(-1)^{k+1}\bE_x^{\pi}\left[Z^k\right]=(-1)^{k+1}\bE_x^{\pi'}\left[Z^k\right], & \textrm{for } k=1,2,\ldots n-1,\\
(-1)^{k+1}\bE_x^{\pi}\left[Z^k\right]>(-1)^{k+1}\bE_x^{\pi'}\left[Z^k\right], & \textrm{for } k=n
\end{cases}
\end{equation}
where $Z:=\sum_{i=0}^{\infty}\beta^i c(X_i,a_i)$ or all moments are equal to each other; for $n=1$ we assume that the first condition in \eqref{eq:moment.optimal} is trivially satisfied. In particular, moment optimal policy $\pi\in\Pi$ automatically satisfies $\pi\in \Pi^*(\beta,0)$ as otherwise condition \eqref{eq:moment.optimal} is directly violated.

\begin{proposition}[A link between moment optimal policy and risk-sensitive optimal policy]\label{pr:moment-optimal}
Fix $\beta\in (0,1)$. Then, the policy $u\in \Pi'$ is moment optimal if and only if there exists $\gamma_\beta<0$ such that $u\in \Pi^{*}(\beta,\gamma)$, for $\gamma\in [\gamma_\beta,0]$.
\end{proposition}

\begin{proof}
The proof follows directly from the power series expansion of the exponent function. For $u,u' \in \Pi'$, $x\in E$, and $\gamma<0$ let
\begin{equation}\label{eq:taylor}
h_x^{u,u'}(\gamma):=\sum_{k=0}^{\infty} \frac{\gamma^k}{k!} \left(\bE_x^{u}\left[Z^k\right]-\bE_x^{u'}\left[Z^k\right]\right).
\end{equation}
note that all moments exist and are finite, so \eqref{eq:taylor} is well posed. Fix $u\in\Pi'$ and assume it is moment optimal. From \eqref{eq:moment.optimal} we get that for any $u'\in\Pi'$ and $x\in E$ the first non-zero term in \eqref{eq:taylor} is even and negative or odd and positive. Consequently, for any $u'\in\Pi'$ and $x\in E$ there exists  $\gamma_{u',x}<0$ such that $h_x^{u,u'}(\gamma)\leq 0$ for $\gamma\in [\gamma_{u',x},0)$; this follow from standard Taylor-based arguments and  Peano reminder analysis. Furthermore, recalling that sets $E$ and $\Pi'$ are finite and setting $\gamma_{\beta}:=\max_{x\in E}\max_{u'\in\Pi'}\gamma_{u',x}<0$, we get that 
\begin{equation}\label{eq:taylor2}
h_x^{u,u'}(\gamma)\leq 0 \quad \textrm{ for any } x\in E,\, u'\in \Pi',\, \gamma\in [\gamma_\beta,0).
\end{equation}
Noting that for $\gamma<0$ the inequality $h_x^{u,u'}(\gamma)\leq 0$ is equivalent to $J_\gamma(x,\pi;\beta)\geq J_\gamma(x,\pi';\beta)$ and that \eqref{eq:moment.optimal} indicates risk-neutral policy optimality (by considering $k=1$), we conclude this part of the proof. The proof of the reverse implication is also based on the analysis of \eqref{eq:taylor}. It is sufficient to notice that \eqref{eq:taylor2} is implying that the first non-zero term in \eqref{eq:taylor} is even and negative or odd and positive, which directly implies \eqref{eq:moment.optimal}.
\end{proof}

In particular, from Proposition~\ref{pr:moment-optimal} we see that by solving the discounted risk-sensitive problem for negative $\gamma$ sufficiently close to zero, we can recover a policy optimal for the risk-neutral discounted problem which admits further properties desirable from the risk management perspective. Of course, similar result could be obtained in the risk-seeking case, by using the standard lexicographic order induced by the moments vector (without -1 multiplication) to define {\it risk-seeking moment optimal} policy.

\section{Numerical examples}\label{S:examples.all}
In this section we present two numerical examples: the first one illustrates the non-stationarity of optimal discounted policies and the dynamics of the turnpike, while the second one shows that the boundary $\gamma(\beta)$ need not be continuous.

\subsection{Independent lottery trials}\label{S:examples}
We illustrate that the optimal policy in the discounted setup could be non-stationary, ultimate stationary policies could be different from Blackwell optimal policies, and investigate the dynamics of the turnpike. We use a simple example based on independent lottery trials based on the idea initially considered in~\cite{Jaq1976} and expanded in \cite{BauPitSte2024,PitSte2025} for a two-action lottery setup. For brevity, we first formulate the initial framework and then consider various parameter setups in specific examples.

Let $(E,U):=(\{1,2,3\},\{a,b,c\})$ and let the transition matrices be given by
\[
P^a=\begin{bmatrix}
0 & 0.2  & 0.8 \\
1 & 0  & 0 \\
1 & 0  & 0 \\
\end{bmatrix},
P^b=\begin{bmatrix}
0 & 0.5  & 0.5 \\
1 & 0 & 0 \\
1 & 0  & 0 \\
\end{bmatrix},
P^c=\begin{bmatrix}
0 & 0.8  & 0.2 \\
1 & 0 & 0 \\
1 & 0  & 0 \\
\end{bmatrix},
\]
for cost function $c\colon E\times U\to\bR$ defined as
\[
c(x,a)=
\begin{cases}
-1& \textrm{if } x=1\\
0& \textrm{if } x=2\\
R& \textrm{if } x=3
\end{cases},\quad
c(x,b)=
\begin{cases}
0& \textrm{if } x=1\\
0& \textrm{if } x=2\\
R& \textrm{if } x=3
\end{cases},\quad
c(x,c)=
\begin{cases}
1& \textrm{if } x=1\\
0& \textrm{if } x=2\\
R& \textrm{if } x=3
\end{cases},
\]
where $R>0$ is a fixed reward; we also fix initial state $X_0=1$. The example could be seen as a sequence of two time-step independent lotteries in which a single trial reward is equal to $R$. For each lottery run, we can choose from three options:
\begin{enumerate}
\item[(a)] we can pay 1 and set the winning chance to 80\%;
\item[(b)] we can have a random setup with 50\% win chance without paying anything;
\item[(c)] we can receive 1 but enter the lottery with the winning chance set to 20\%.
\end{enumerate}
Consequently, every second step we need to make a decision which lottery option to choose. It should be noted that while the example does not satisfy the global mixing assumptions \eqref{A.1} and \eqref{A.2}, the realised multi-step transition probabilities are strongly mixing if we restrict ourselves to the MDPs with the initial state $X_0=1$. Furthermore, one can easily modify the example without materially impacting the results by slightly modifying the transition probabilities to ensure that all matrix entries are positive, see \cite{PitSte2025} for details.

Let us now provide a more formal comment. In general, we have $27$ Markov stationary policies. Noting that for $i\in\{2,3\}$ we have $\bP^a(i,\cdot)\equiv \bP^b(i,\cdot)\equiv \bP^c(i,\cdot)$ and  $c(i,\cdot)\equiv \textrm{const}$, it does not effectively matter which action is assigned to state $2$ and $3$.  Consequently, it is sufficient to fix three stationary Markov policies $u_a,u_b,u_c\in \Pi'$ such that $u_a(\cdot)\equiv a$, $u_b(\cdot)\equiv b$, and $u_c(\cdot)\equiv c$, and consider Markov policies of the form $\pi=(a_i)_{i\in\bN}$, where $a_{2k}\in \{u_a, u_b, u_c\}$, and $a_{2k+1}\equiv a_{2k}$, for $k\in\bN$. In that case, recalling that entropic operator is additive with respect to independent random variables and each lottery choice could be made independently of the past choices (as the state resets to 1 every second step), we get
\begin{align}
J_{\gamma}(1,\pi;\beta) & =\tfrac{1}{\gamma}\ln\bE_1^{\pi}\left[\exp\left(\gamma\sum_{i=0}^{\infty}\beta^{i} c(X_i,a_i)\right)\right]\nonumber\\
&=\tfrac{1}{\gamma}\ln\bE_1^{\pi}\left[\exp\left(\gamma\sum_{k=0}^{\infty}\beta^{2k} (c(X_{2k},a_{2k})+\beta c(X_{2k+1},a_{2k+1}))\right)\right] \nonumber\\
&=\tfrac{1}{\gamma}\ln\bE_1^{\pi}\left[\exp\left(\gamma\sum_{k=0}^{\infty}\beta^{2k} ([\1_{\{a_{2k}=u_c\}}-\1_{\{a_{2k}=u_a\}}]+\beta R\cdot\1_{\{X_{2k+1}=3\}})\right)\right]\label{eq:lottery1}
\end{align}
Noting that the constituents of the sum in \eqref{eq:lottery1} form a sequence of affine transformed independent Bernoulli random variables, we can rewrite \eqref{eq:lottery1} as 
\begin{equation}\label{eq:lottery2}
J_{\gamma}(1,\pi;\beta)=\sum_{k=0}^{\infty}\left(\beta^{2k}c_k+\tfrac{1}{\gamma}\ln\left[p_k e^{0} +(1-p_k)e^{\gamma\beta^{2k+1} R}\right]\right),
\end{equation}
where $c_k:=\1_{\{a_{2k}=u_c\}}-\1_{\{a_{2k}=u_a\}}$ and $p_k:=0.5 +0.3c_k$. Consequently, to find the maximum of \eqref{eq:lottery2} it is sufficient to (independently) maximise each constituent of the sum in \eqref{eq:lottery2}, that is, the value
\begin{equation}\label{eq:lottery.Ak}
A_k(u,\gamma,\beta):=\left(\beta^{2k}c_u+\tfrac{1}{\gamma}\ln\left[p_u e^{0} +(1-p_u)e^{\gamma \beta^{2k+1} R}\right]\right), \textrm{ for } u\in\{u_a,u_b,u_c\},
\end{equation}
where $c_u:=\1_{\{u=c\}}-\1_{\{u=a\}}$, $p_u:=0.5 +0.3c_u$, and $k\in\bN$. Note that the three possible choices correspond to lottery policies $a$, $b$, and $c$ as explained in the non-formal example's introduction. Using similar calculations, one can also show that for the risk-sensitive averaged criterion \eqref{eq:RSC.averaged}, for Markov stationary policies, we get
\begin{equation}\label{eq:J.averaged.lottery}
\tilde J_\gamma(1,u)=\frac{1}{2}\left( c_u+\tfrac{1}{\gamma}\ln\left[p_u e^{0} +(1-p_u)e^{\gamma R}\right]\right),
\end{equation}
where $c_u:=\1_{\{u=c\}}-\1_{\{u=a\}}$ and $p_u:=0.5 +0.3c_u$. For stationary Markov policies, one can also calculate the value for the discounted risk-neutral framework (for $\gamma=0$) and get
\begin{equation}\label{eq:J.neutral.lottery}
J_0(1,u;\beta)=\frac{1}{1-\beta^2}(c_u +(1-p_u)\beta R).
\end{equation}
We are now ready to consider specific parameter choices that will illustrate the concepts considered in this paper.

\begin{example}[Optimal non-stationary policy in the discounted framework]\label{ex:1}
 Let us fix the reward $R=7$, risk sensitivity $\gamma=-1$ and discount factor $\beta=0.95$. Then, we can compute the values of $A_k(\cdot,\gamma,\beta)$ to determine the optimal discounted policy. The values for $k=1,2\ldots,10$ are presented in Table~\ref{T:ex1}. As can be seen, the optimal discounted policy is not stationary and it is best to choose policy $c$ in the first four runs and then switch to policy $a$. By direct computations, we can check that for the policy $\pi$ given by
\[
\pi=(\underbrace{u_c,u_c,u_c,u_c,u_c,u_c,u_c,u_c}_\text{4 lottery runs (8 steps)},u_a,u_a,u_a,u_a,\ldots),
\]
we get $J_{\gamma}(1,\pi;\beta)\approx 22.8$ while for stationary Markov policies we get $J_{\gamma}(1,u_a;\beta)\approx 21.4$,  $J_{\gamma}(1,u_b;\beta)\approx 15.5$, and $J_{\gamma}(1,u_c;\beta)\approx 15.6$, so that this non-stationary Markov policy is superior when compared to all stationary policies. By direct comparison of $A_k(\cdot,\gamma,\beta)$ for large $k$s, one can show that the policy $\pi$ is in fact optimal, see \cite{Jaq1976} for the analytical-based argumentation idea.

\begin{table}[ht]
\centering
{\small
\begin{tabular}{cr|rrrrrrrrrr}
  \hline
$k$ &  & 0 & 1 & 2 & 3 & 4 & 5 & 6 & 7 & 8 & 9 \\ 
  \hline
\multirow{3}{*}{$A_k$} &a & 0.60 & 0.70 & 0.78 & 0.84 & {\bf 0.90} & {\bf 0.94} & {\bf 0.96} & {\bf 0.98} & {\bf 0.98} & {\bf 0.96} \\ 
 & b & 0.69 & 0.69 & 0.69 & 0.69 & 0.68 & 0.67 & 0.67 & 0.65 & 0.64 & 0.62 \\ 
  &c & {\bf 1.22} & {\bf 1.13}  & {\bf 1.04} & {\bf 0.96} & 0.88 & 0.82 & 0.76 & 0.70 & 0.65 & 0.60 \\   \hline
  \multicolumn{2}{l|}{Best}   & c & c & c & c & a & a & a & a & a & a \\ 
   \hline
\end{tabular}}
\caption{The table presents the first ten values of $A_k(u,\gamma,\beta)$, for policy choices $u\in\{u_a,u_b,u_c\}$ under the setup from Example~\ref{ex:1} with $R=7$, $\beta=0.95$ and $\gamma=-1$. As one can see, the best choice is to use policy corresponding to lottery $c$ for the first four lotteries, and then switch to lottery $a$. This shows that the optimal policy for the corresponding discounted risk-sensitive objective criterion must be a non-stationary Markov policy.}\label{T:ex1}
\end{table} 
\end{example}

\begin{example}[Difference between ultimately stationarity and Blackwell optimality]\label{ex:2}
Let us fix the reward $R=7$ and numerically analyze the turnpike values $N(\beta,\gamma)$, that is, the earliest ultimate stationarity threshold. We consider the compact parameter set for $\gamma\in [-2.5,2.5]$ and $\beta\in [0.9,0.995]$ and calculate the values of $N(\beta,\gamma)$ on a $(0.001)$-step parameter grid by comparing the values of $A_k(u,\gamma,\beta)$ and noting when the last switch takes place; see \cite{PitSte2025} for a similar level-plot setup where more details are provided. The level plot graph is presented in Figure~\ref{F:ex2}.

\begin{figure}[htp!]
  \centering
  \fbox{%
    \includegraphics[width=0.45\linewidth]{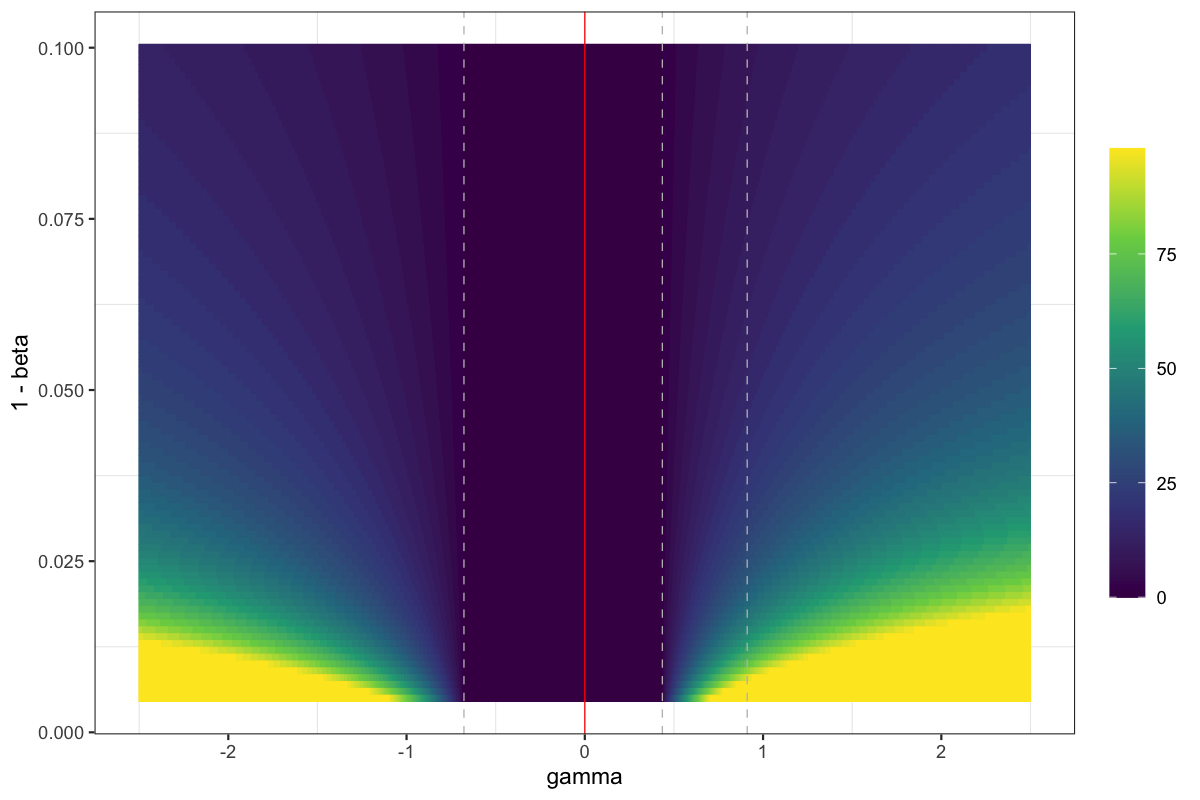}%
  }
  \caption{The plot presents the (ultimately stationary) optimal policy turnpike $N(\beta,\gamma)$ for different values of $\gamma\in[-2.5,2.5]$ and $\beta\in(0.9,0.995)$ under the setup in Example~\ref{ex:2}, for $R=7$. The vertical dashed lines mark $\gamma$-based thresholds when the averaged per unit of time optimal (stationary) policy is changing. As one can see, we can observe that when $\gamma\to 0$, then the turnpike approaches zero leading to stationary optimal policy. On the other hand, when $\beta\to 1$, the turnpike could approach either 0 or $+\infty$ depending on the limit policy for the averaged problem.}
  \label{F:ex2}
\end{figure}

From the plot we can deduce that for small values of $\gamma$ one can find a policy that is globally stationary for all values of $\beta\in [0.9,0.995]$. In fact, one can show that for $\gamma\in [-0.67,0.43]$ the ultimate stationary policy is equal to the averaged risk sensitive stationary policy. This can be illustrated by analyzing the behavior of the problem limit when $\beta\to 1$ or $\gamma\to 0$ and computing the corresponding optimal values \eqref{eq:J.averaged.lottery} and \eqref{eq:J.neutral.lottery} for all stationary policies. The illustration is presented in Figure~\ref{F:ex2b}. 

\begin{figure}[htp!]
  \centering
  \fbox{%
    \includegraphics[width=0.45\linewidth]{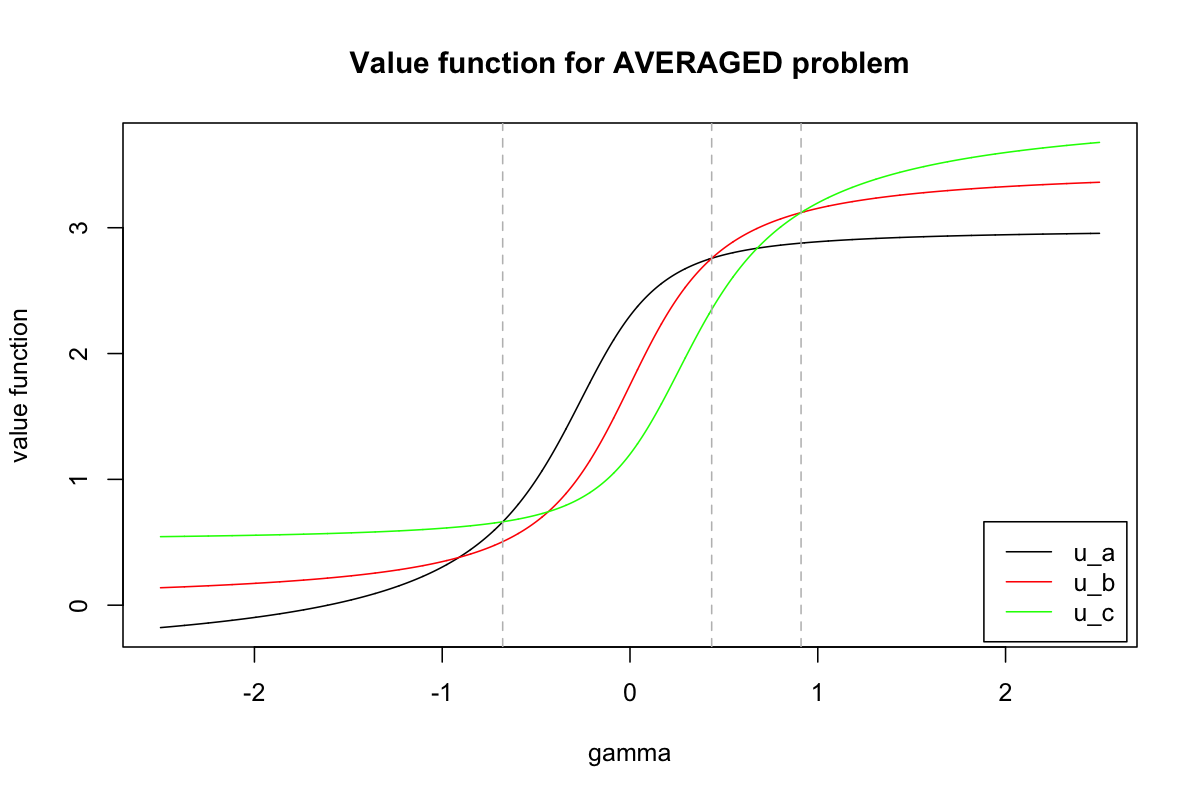}%
     \includegraphics[width=0.45\linewidth]{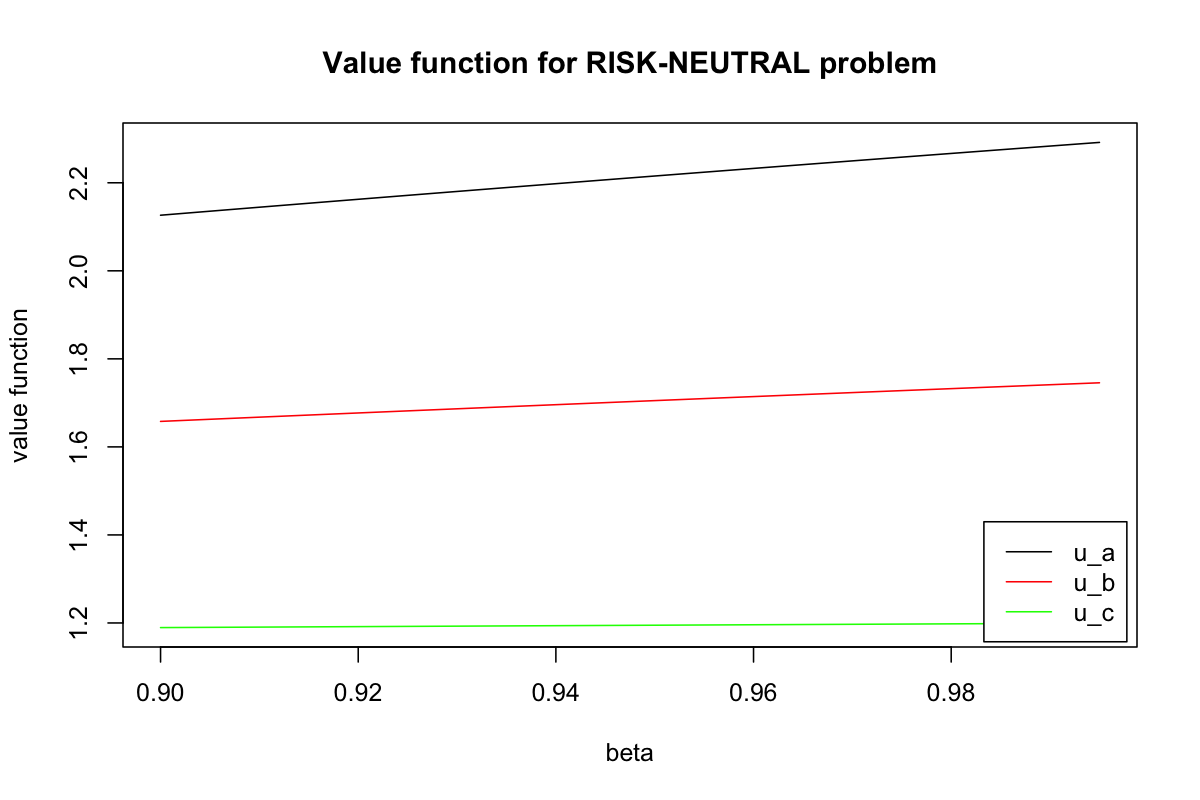}%
  }
  \caption{The plot presents the value function for the stationary policies in the limit case when $\beta\to 1$ (risk-averaged) and when $\gamma\to 0$ (risk-neutral) using the value functions \eqref{eq:J.averaged.lottery} and \eqref{eq:J.neutral.lottery}; we have multiplied the second function by $(1-\beta)$ to normalise the output as otherwise it tends to $+\infty$ as $\beta\to 1$. The results are presented for different values of $\gamma\in[-2.5,2.5]$ and $\beta\in(0.9,0.995)$ under the setup in Example~\ref{ex:2}, for $R=7$. The vertical dashed lines mark $\gamma$-based thresholds when the averaged per unit of time optimal (stationary) policy is changing. Note that in the considered set, the policy $u_a$ is risk-neutral discounted optimal for all values of $\beta\in [0.9,0.995]$.}
  \label{F:ex2b}
\end{figure}

As expected, for $\gamma\in [-0.67,0.43]$, the policy $u_a$ is indeed optimal for the averaged problem, and for $\beta\in [0.9,0.995]$, $u_a$ is also optimal for the risk-neutral problem, which results in a stable stationary policy region. On the other hand, when $\gamma\not\in [-0.67,0.43]$, the Blackwell-optimal stationary policy is not equal to the tail stationary policy coming from the ultimately stationary optimal policy representation. For more detailed numerical analysis and comment on the structural differences between ultimately stationary policies and Blackwell optimal policies, we refer to \cite{BauPitSte2024} and \cite{PitSte2025}.
\end{example}

\begin{example}[Turnpike dynamics]\label{ex:3}
In the setup considered in Example~\ref{ex:2} the risk-neutral discounted value function linked to $u_a$ was dominant, leading to relatively stable turnpike behavior; see Figure~\ref{F:ex2b} for the dynamics of the value function and Figure~\ref{F:ex2} for the turnpike dynamics. To see how sensitive the dynamics of the turnpike could be, let us decrease the reward value  and consider $R=3.5$. In this setup, one can show that the optimal stationary policy for the discounted risk-neutral problem is switching on the considered parameter grid, that is, when $\beta\in [0.9,0.995]$. To illustrate this, we present the analog of the plot in Figure~\ref{F:ex2b} for the adjusted setup, see Figure~\ref{F:ex3b}.

\begin{figure}[htp!]
  \centering
  \fbox{%
    \includegraphics[width=0.45\linewidth]{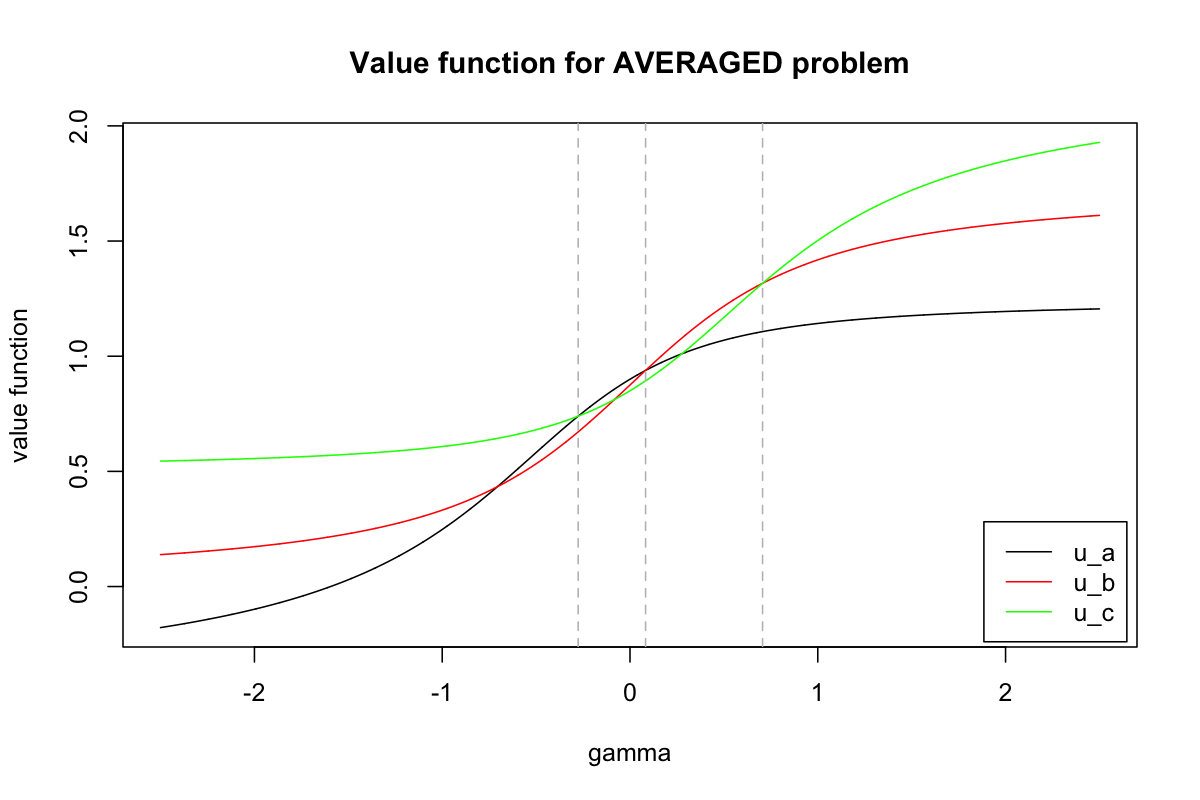}%
     \includegraphics[width=0.45\linewidth]{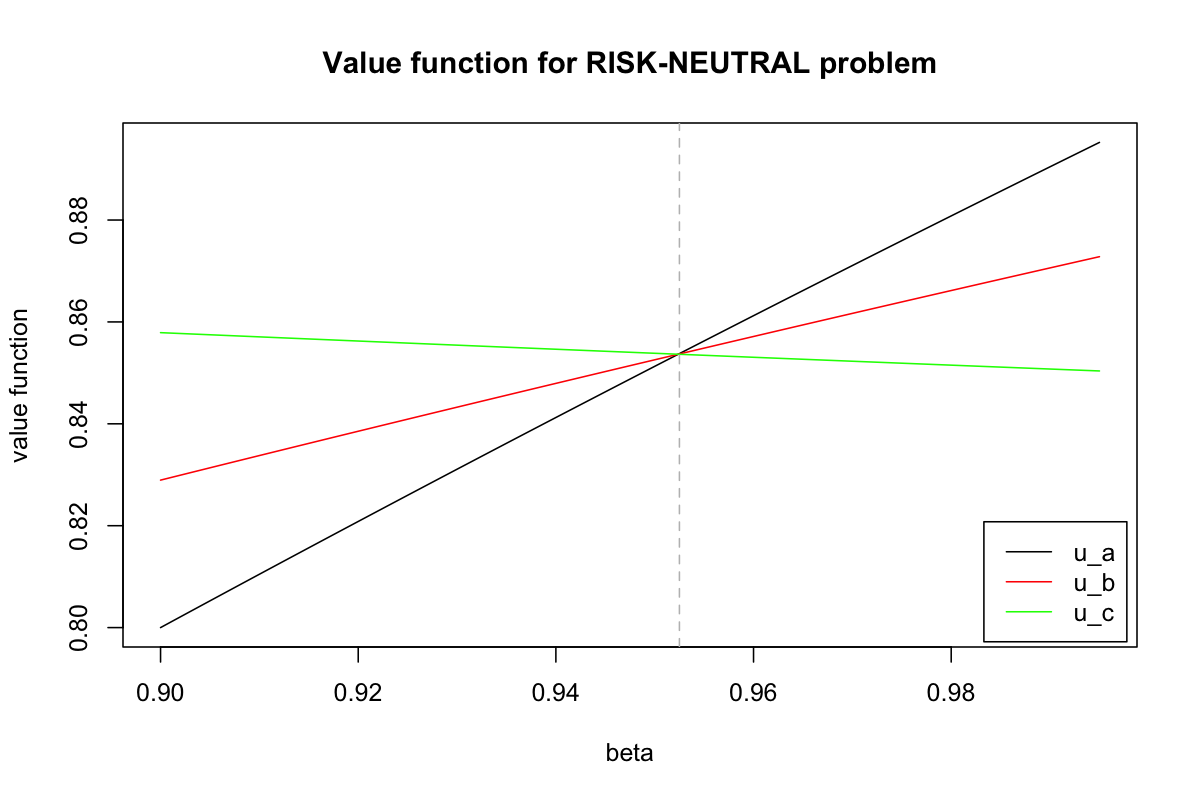}%
  }
  \caption{The plot presents the value function for the stationary policies in the limit case when $\beta\to 1$ (risk-averaged) and when $\gamma\to 0$ (risk-neutral) for the decreased reward $R=3.5$ in the Example~\ref{ex:2} framework; see Figure~\ref{F:ex2b} caption for detailed description. Note that the risk-neutral optimal discounted policy is changing - the switch is at $\beta\approx 0.9525$.}
  \label{F:ex3b}
\end{figure}

As one can see, the risk-neutral discounted policy is indeed switching for $\beta\approx 0.9525$ and the value functions are close to each other in the neighborhood. To see how this impacts the dynamics of the turnpike, let us plot the analogue of Figure~\ref{F:ex2}.

\begin{figure}[htp!]
  \centering
  \fbox{%
    \includegraphics[width=0.45\linewidth]{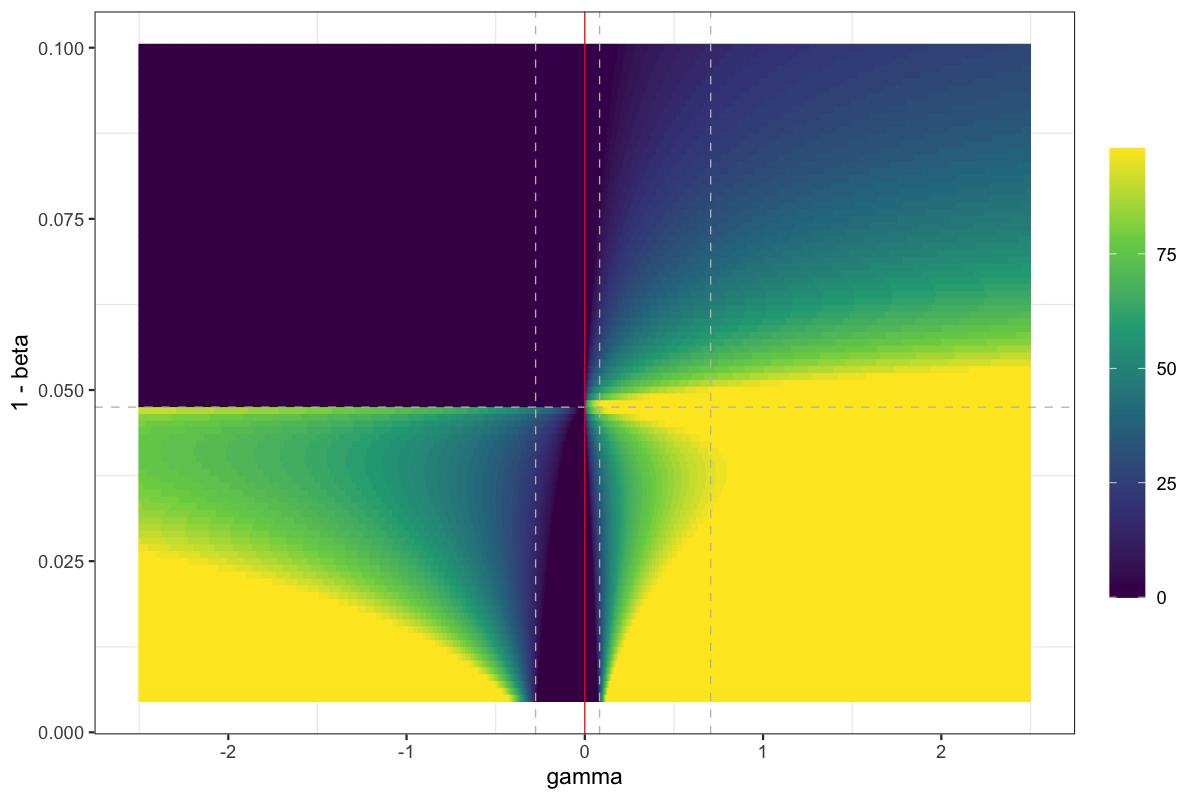}%
  }
  \caption{The plot presents the (ultimately stationary) optimal policy turnpike $N(\beta,\gamma)$ for different values of $\gamma\in[-2.5,2.5]$ and $\beta\in(0.9,0.995)$ under the setup in Example~\ref{ex:3}, with decreased reward $R=3.5$; see Figure~\ref{F:ex2} caption for detailed description. Horizontal dashed line indicate optimal stationary risk-neutral discounted policy switch point.  
  }
  \label{F:ex3}
\end{figure}

The result is presented in Figure~\ref{F:ex3}. As we can see, the dynamics of the turnpike is materially different - we observe that the value of the turnpike need not be monotone with respect to $\beta$ and could exhibit a materially different behavior depending on the interaction between the limit policies.
\end{example}

\subsection{A single lottery choice and discontinuity of the boundary $\gamma(\beta)$}\label{sec:Exgamma}

In this section we illustrate the behaviour of the boundary $\gamma(\beta)$ defined in \eqref{eq:gamma.beta} and show that the mapping $\beta\mapsto\gamma(\beta)$ need not be continuous, so that the lower semicontinuity stated in Corollary~\ref{cor:boundary}~b) cannot be improved. We follow a single lottery choice setup, using logic similar to the one from Section~\ref{S:examples}.

Let $(E,U):=(\{1,\ldots,7\},\{a,b\})$ and let the transition probabilities be given by $\bP^{a}(1,\cdot)=\big(\tfrac{1}{2},\tfrac{1}{2}\big)$ on $\{2,3\}$ and $\bP^{b}(1,\cdot)=\big(\tfrac{1}{10},\tfrac{9}{10}\big)$ on $\{5,6\}$, while the transitions in the remaining states are deterministic and do not depend on the action; namely, we fix $2,3\mapsto 4$, $4\mapsto 7$, $5,6\mapsto 7$, and $7\mapsto 7$. The reward function $c\colon E\times U\to\bR$ is action-dependent only in the state $1$ and given by
\[
c(1,a)=\tfrac{1}{4},\quad c(1,b)=0,\quad
c(2,\cdot)=-2,\quad c(3,\cdot)=0,\quad c(4,\cdot)=1,
\]
\[
c(5,\cdot)=-\tfrac{9}{5},\quad c(6,\cdot)=\tfrac{1}{5},\quad c(7,\cdot)=0;
\]
we also fix the initial state $X_0=1$. The example can be seen as a single choice between two fair lotteries. Namely, let $\xi_a$ and $\xi_b$ be random variables with
\begin{equation}\label{eq:ex4.xi}
\bP(\xi_a=-1)=\bP(\xi_a=1)=\tfrac{1}{2},\qquad
\bP(\xi_b=-\tfrac95)=\tfrac{1}{10},\quad \bP(\xi_b=\tfrac15)=\tfrac{9}{10},
\end{equation}
so that $\bE[\xi_a]=\bE[\xi_b]=0$, $\textnormal{Var}(\xi_a)=1$, and $\textnormal{Var}(\xi_b)=0.36$; the second lottery is thus less volatile than the first one, but it admits a heavier loss. For the two available actions we can:
\begin{enumerate}
\item[(a)] we receive $\tfrac14$ immediately, then in the next period we pay $1$ and collect the
payoff $\xi_a$ of the first lottery, and one period later we receive $1$;
\item[(b)] in the next period we collect the payoff $\xi_b$ of the second lottery, at no cost.
\end{enumerate}
Since the state $1$ is visited only at time $0$ and the process is absorbed in the state $7$ afterwards, every policy $\pi\in\Pi$ is, from the state $1$, equivalent to one of the two stationary policies $u_a,u_b\in\Pi'$ which choose $a$, respectively $b$, in the initial state $1$; in particular, an optimal stationary policy exists for every $(\beta,\gamma)$, and the one-step deviation differences from Section~\ref{S:small.aversion} reduce to the comparison of $u_a$ and $u_b$. Note that, as in Section~\ref{S:examples}, the example does not satisfy \eqref{A.1} and \eqref{A.2}, which are not used in Section~\ref{S:small.aversion}. By simple calculations, for $\gamma\neq0$, we get
\begin{align*}
J_{\gamma}(1,u_a;\beta) &=\tfrac1\gamma \ln \bE\big[\exp\big(\gamma\big(\tfrac14+\beta(-1+\xi_a)+\beta^2\big)\big)\big],\\
J_{\gamma}(1,u_b;\beta) &=\tfrac1\gamma \ln \bE\big[\exp\big(\gamma\beta\xi_b\big)\big].
\end{align*}
Recalling that the entropic utility $\mu^{\gamma}$ is cash-additive (see Appendix~\ref{A:entropic}), we get
\begin{equation}\label{eq:ex4.values}
J_{\gamma}(1,u_a;\beta)=\big(\beta-\tfrac12\big)^2+\beta\mu^{\beta\gamma}(\xi_a),\qquad
J_{\gamma}(1,u_b;\beta)=\beta\mu^{\beta\gamma}(\xi_b),
\end{equation}
where $\beta\mu^{\beta\gamma}(\xi_a)=\frac{1}{\gamma}\ln\cosh(\beta\gamma)$ and $\beta\mu^{\beta\gamma}(\xi_b)=\frac{\beta}{5}+\frac{1}{\gamma}\ln\left(\frac{9+e^{-2\beta\gamma}}{10}\right)$.

\begin{example}[Discontinuity of the boundary $\gamma(\beta)$]\label{ex:4}
By \eqref{eq:ex4.values}, for $\gamma<0$, we get $u_a\in\Pi^{*}(\beta,\gamma)$ if and only if
\begin{equation}\label{eq:ex4.constraints}
\beta\left[\mu^{\beta\gamma}(\xi_b)-\mu^{\beta\gamma}(\xi_a)\right]\leq \big(\beta-\tfrac12\big)^2 .
\end{equation}
Recalling that $\bE[\xi_a]=\bE[\xi_b]=0$ and using the expansion $\mu^{t}(\xi)=\tfrac{t}{2}\textnormal{Var}(\xi)+O(t^2)$, we get that the left-hand side of \eqref{eq:ex4.constraints} equals $\tfrac{\beta^2|\gamma|}{2}\left(\textnormal{Var}(\xi_a)-\textnormal{Var}(\xi_b)\right)(1+o(1))$ as $\gamma\uparrow 0$; in particular, it is strictly positive for $\gamma<0$ sufficiently close to zero, as the first lottery is more volatile than the second one. Moreover, in the risk-neutral case ($\gamma=0$) we get $J_0(1,u_a;\beta)=(\beta-\tfrac12)^2$ and $J_0(1,u_b;\beta)=0$, so that $u_a$ is the unique risk-neutral optimal policy for $\beta\neq\tfrac12$, whereas for $\beta=\tfrac12$ both policies are optimal. Let us now study continuity of $\beta\mapsto\gamma(\beta)$ at point $\beta=\tfrac{1}{2}$.

Since $u_b\notin\Pi^{*}(\beta,0)$ for $\beta\neq\tfrac12$, we get $\rho_{u_b}(\beta)=+\infty$ and consequently $\gamma(\beta)=\rho_{u_a}(\beta)$, that is, $\gamma(\beta)$ is the smallest $\gamma\leq0$ for which \eqref{eq:ex4.constraints} holds on the whole interval $[\gamma,0]$. In particular, $0\geq\gamma(\beta)\geq\gamma_1(\beta)$, where $\gamma_1(\beta)$ denotes the largest negative solution of
\begin{equation}\label{eq:ex4.curve}
\beta\left[\mu^{\beta\gamma}(\xi_b)-\mu^{\beta\gamma}(\xi_a)\right]= \big(\beta-\tfrac12\big)^2 ,
\end{equation}
which exists for $\beta$ sufficiently close to $\tfrac12$. By direct calculations we get $\gamma_1(\beta)=-\tfrac{25}{2}(\beta-\tfrac12)^2(1+o(1))$ as $\beta\to\tfrac12$, which implies
\begin{equation}\label{eq:ex4.limit}
\lim_{\beta\to 1/2}\gamma(\beta)=0.
\end{equation}
Now, for $\beta=\tfrac12$, the right-hand side of \eqref{eq:ex4.constraints} vanishes, while its left-hand side is strictly positive for $\gamma<0$ sufficiently close to zero. Consequently, $u_a\notin\Pi^{*}(\tfrac12,\gamma)$ for such $\gamma$ and $u_b$ is the unique optimal policy on a non-degenerate interval of risk-aversion parameters. Recalling that $u_b$ is also risk-neutral optimal for $\beta=\tfrac12$, we get
\begin{equation}\label{eq:ex4.jump}
\gamma(\tfrac12)\leq\rho_{u_b}(\tfrac12)<0,
\end{equation}
which shows that $\gamma(\tfrac{1}{2})\neq 0$; the heavier loss of the second lottery makes $\rho_{u_b}(\tfrac12)$ finite, and numerically we get $\gamma(\tfrac12)\approx-3.49$. Combining \eqref{eq:ex4.limit} with \eqref{eq:ex4.jump} we conclude that $\lim_{\beta\to1/2}\gamma(\beta)>\gamma(\tfrac12)$, that is, the mapping $\beta\mapsto\gamma(\beta)$ is lower semicontinuous, in accordance with Corollary~\ref{cor:boundary}~b), but it is not continuous at $\beta=\tfrac12$; see Figure~\ref{F:ex4} for an illustration. Let us also note that \eqref{eq:ex4.limit} shows that the constant $\underline\rho$ vanishes on every parameter rectangle containing $\beta=\tfrac12$, so that the assumption $\underline\rho>0$ made in Proposition~\ref{pr:finite.partition} is not automatically satisfied.

\begin{figure}[htp!]
  \centering
  \fbox{%
    \includegraphics[width=0.45\linewidth]{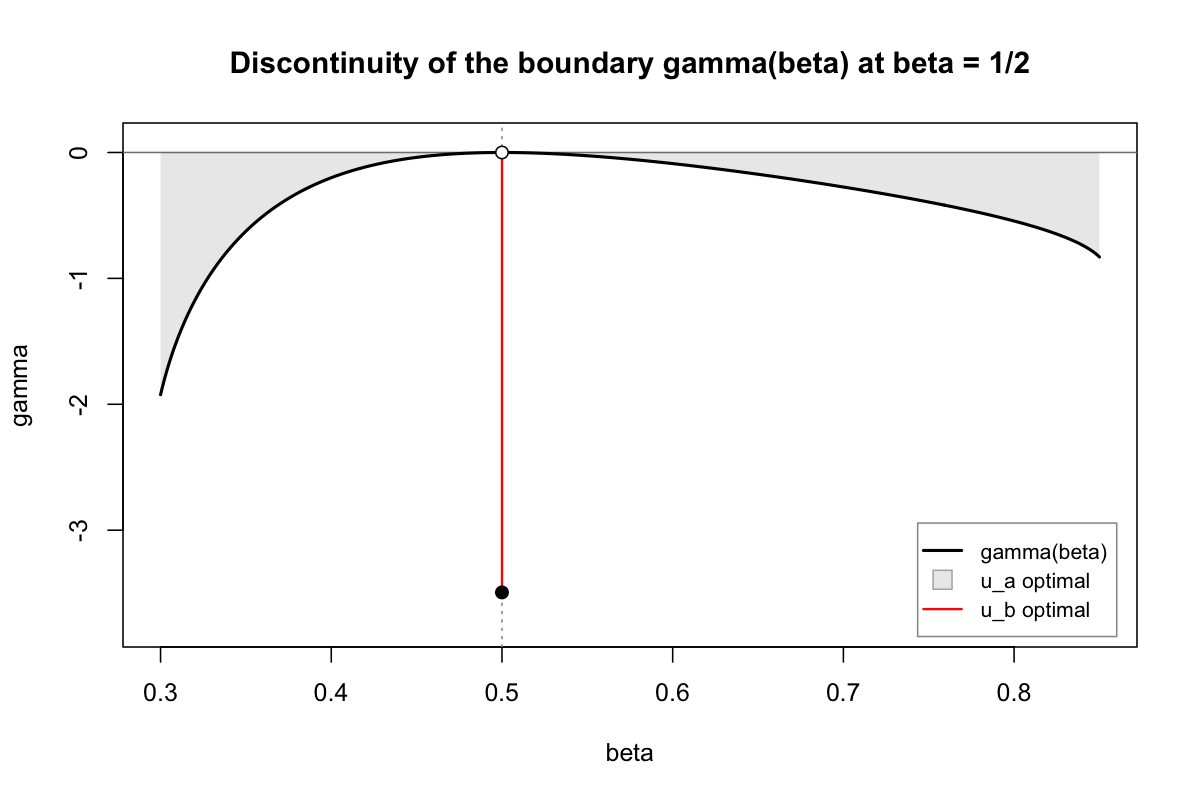}%
  }
  \caption{The plot presents   $\beta\mapsto \gamma(\beta)$  under the setup in Example~\ref{ex:4}. In the shaded region the policy $u_a$ is optimal, while at $\beta=\tfrac12$ the policy $u_b$ is optimal on the marked vertical segment. As one can see, the boundary collapses to zero as $\beta\to\tfrac12$, while at $\beta=\tfrac12$ it drops to $\gamma(\tfrac12)\approx-3.49$; at $(\beta,\gamma)=(\tfrac12,0)$ both policies are optimal.}
  \label{F:ex4}
\end{figure}
\end{example}

\appendix

\section{Selected properties of entropic utility measure}\label{A:entropic}
In this section, we recall some basic properties of  entropic utility measures. For brevity, we use $(\Omega,\cF,\bP)$ to denote a generic probability space with the expectation operator $\bE$. For any $\gamma\in\bR$ and $X\in L^{\infty}(\Omega,\cF,\bP)$ we use
\[
\mu^{\gamma}(X):=
\begin{cases}
\tfrac{1}{\gamma}\ln \bE[e^{\gamma X}], & \textrm{if } \gamma\neq 0,\\
\bE[X], &\textrm{if } \gamma=0,
\end{cases}
\]
to denote {\it entropic utility measure} with risk-aversion $\gamma$.

\begin{proposition}\label{pr:entropic.monotone}
Let $X,Y\in L^{\infty}(\Omega,\cF,\bP)$. Then

\begin{enumerate}

\item[a)] If $X\leq Y$, then $\mu^{\gamma}(X)\leq \mu^{\gamma}(Y)$, for any $\gamma\in\bR$.
\item[b)] We have $-\|X\|_{\infty} \leq \mu^{\gamma}(X) \leq \|X\|_{\infty}$, for $\gamma\in\bR$.
\item[c)] For any $\gamma\neq 0$, we get $\mu^{\gamma}(X)\leq \bE[X]+\tfrac{|\gamma|}{2}\|X\|_{\textnormal{sp}}^2$.
\item[d)] The mapping $\gamma\to \mu^{\gamma}(X)$ is monotone and real-analytic.
\end{enumerate}
\end{proposition}

\begin{proof}

For brevity, we omit the proof of properties a), b). Property c) follows directly from Hoeffding's Lemma (see \cite{BouLugMas2013}), that is, inequality
\begin{equation}\label{eq:hoef}
\bE[e^{\gamma X}] \leq e^{\gamma \bE[X] +\tfrac{1}{8}\gamma^2(\esssup X-\essinf X)^2},
\end{equation}
for any $\gamma\in \bR$. Indeed, it is sufficient to note that for $\gamma<0$ we get $\mu^\gamma(X)\leq \bE[X]$ from the monotonicity of entropic utility, while for $\gamma>0$, \eqref{eq:hoef} implies directly $\mu^{\gamma}(X)\leq \bE[X]+\tfrac{\gamma}{2}\|X\|_{\textnormal{sp}}^2$. For brevity, we also omit the proof of monotonicity of $\gamma\to \mu^{\gamma}(X)$ in d); see \cite{KupSch2009}. To show that $\gamma\to \mu^{\gamma}(X)$ is real-analytic, we first observe that $\gamma\to \bE[e^{\gamma X}]$ is a moment generating function of a bounded random variable which is real-analytic, see \cite{Cur1942}. The function remains real-analytic after taking the logarithm; see Proposition 1.4.2 in \cite{KraPar2002}. Multiplication by $\tfrac{1}{\gamma}$ also preserves real-analytic property on $(-\infty,0)$ and $(0,+\infty)$. Consequently, it is sufficient to show that $\gamma\to \tfrac{1}{\gamma}\ln\bE[e^{\gamma X}]$ is real-analytic in 0; recall that we have already shown that the function is continuous. As $\gamma\to \ln\bE[e^{\gamma X}]$ is analytic in $\bR$ and $\ln\bE[e^{0}]=0$, we get that $a_0=0$ in local expansion $\ln\bE[e^{\gamma X}]=\sum_{i=0}^{\infty} a_i \gamma^i$. This implies that $\tfrac{1}{\gamma}\ln\bE[e^{\gamma X}]=\sum_{i=0}^{\infty} a_{i+1} \gamma^{i}$ is well defined, which concludes this part of the proof.
\end{proof}

\section*{Declaration}
The authors used large language model assistants (Anthropic Claude Opus~5; OpenAI ChatGPT 5) in July 2026 to generate the example in Section~\ref{sec:Exgamma}. All statements and calculations have been independently verified and proved by the authors, who assume responsibility for all content.

\bibliographystyle{siamplain}
\bibliography{references}
\end{document}